\documentclass[a4paper,11pt]{article}

\usepackage{latexsym}
\usepackage{amsmath,amsfonts,amssymb,amsthm}
\newtheorem{thm}{Theorem}[section]

\newtheorem{cor}[thm]{Corollary}
\newtheorem{prop}[thm]{Proposition}

\newtheorem{rem}[thm]{Remark}

\input epsf

\title{Determinants related to Dirichlet characters modulo $2,\ 4$ and
$8$ of binomial coefficients and the algebra of recurrence matrices}
\author{Roland Bacher}
\date{}
\begin{document}
\maketitle
%\par qbin10.tex dans recherche/thue-morse

{\sl Abstract: Using recurrence matrices, defined and described with some
details, we study a few determinants related to evaluations of
binomial coefficients on Dirichlet characters modulo $2,4$ and $8$.}
\footnote{Math. Class: 
05A10,05A30,11B65,11B85,15-99 Keywords: 
binomial coefficient, $q-$binomial coefficient, Dirichlet character,
determinant, recurrence matrix} 

\section{Introduction}

The aim of this paper is twofold: It
contains computations of a few determinants
related to binomial coefficients. The most interesting example, discovered 
after browsing through \cite{Gr1}, is obtained by considering 
binomial coefficients modulo $4$.
We include two similar examples taken from \cite{BacherChapman}
and \cite{BacherCRAS}.

The second topic discussed in this work are recurrence matrices,
see \cite{BacherCRAS} for a very condensed outline. They are defined
as certain sequences of matrices involving self-similar structures
and they form an algebra. Computations in the algebra of recurrence 
matrices are the main tool for proving the determinant formulae mentionned 
above. Recurrence matrices are however of independent interest 
since they are closely
linked for example to automatic sequences, see \cite{AS}, to rational formal
power series in free noncommuting variables and to
groups of automata, see \cite{Grigorchuk} for the perhaps most 
important example.
The determinant calculations of this paper can thus be considered as
illustrations of some interesting features
displayed by recurrence matrices.

The sequel of this paper is organised as follows:

The next section recalls mostly well-known facts concerning 
binomial and $q-$binomial coefficients and states the main results.

Section \ref{sectalgrec} defines the algebra $\mathcal R$ of 
recurrence matrices. It describes them with more details than 
necessary for proving the formulae of Section \ref{sectMain}.

Section \ref{sectinv} discusses a few features of the
group of invertible elements in $\mathcal R$.

Section \ref{sectmod2} proves formulae for the determinant of the
reduction modulo $2$ of the
symmetric Pascal matrix (already contained in \cite{BacherChapman})
and of a determinant related to the $2-$valuation of the binomial
coefficients (essentially contained in 
\cite{BacherCRAS}).

Section \ref{sectBeeblebrox} is devoted to the proof of our main
result, a formula for $\det(Z(n))$ where $Z(n)$
is the matrix with coefficients $\chi_B({s+t\choose s})\in\{0,\pm 1\},
\ 0\leq s,t<n$
obtained by considering the ``Beeblebrox reduction''
(given by the Dirichlet character $\chi_B(2m)=0,\ 
\chi_B(4m\pm 1)=\pm 1$ for $m\in\mathbb Z$) of binomial coefficients.
The proof uses an $LU$ factorisation of the infinite symmetric
matrix $Z=Z(\infty)$ 
and suggests to consider two (perhaps interesting) groups 
$\Gamma_L$ and $\Gamma_Z$ whose generators display beautiful
``self-similar'' structures.
This section ends with a short digression on the ``lower triangular
Beeblebrox matrix'' and the associated group.

Section \ref{sectjacDir} contains some data concerning
the reduction of binomial coefficients by a Dirichlet 
character modulo $8$ related to the Jacobi symbol.

Section \ref{sectqbin} reproves a known formula
for evaluating $q-$binomial coefficients at roots of unity. 
This formula yields easily formulae for some determinants associated to 
the reduction modulo $2$ and the Beeblebrox reduction of (real
and imaginary parts) of $q-$binomial coefficients evaluated at $q=-1$ and
$q=i$.

\section{Main results}\label{sectMain}

\subsection{Reductions modulo 2}\label{subsectmod2}

Let $P(n)$ be the integral symmetric $n\times n$ matrix with coefficients
$P_{s,t}\in\{0,1\},\ 0\leq s,t<n$ defined by
$$P_{s,t}\equiv {s+t\choose s}\pmod 2$$
where ${s+t\choose s}=\frac{(s+t)!}{s!\ t!}$ denotes the usual binomial
coefficient involved in the expansion $(x+y)^n=\sum_{k=0}^n
{n\choose k}x^ky^{n-k}$.

The evaluation ${s+t\choose s}\pmod 2$ can be computed 
using a Theorem of Lucas, see \cite{Lucas}, page 52. Given 
a prime number $p$, it states that
$${n\choose k}\equiv \prod_{j\geq 0}{\nu_j\choose \kappa_j}\pmod p$$
where $\nu_i,\kappa_i\in\{0,1,\dots,p-1\}$ are the coefficients 
of the $p-$ary expansion of $n=\sum_{j\geq 0}\nu_j p^j$ and 
$k=\sum_{j\geq 0}\kappa_j p^j$. Another formula
(due to Kummer) for ${n\choose k}\pmod 2$ will be presented in 
Section \ref{subsec2val}.

Let $\mathop{ds}(n)=\sum_{j=0}^{\lfloor \mathop{log}_2(n)\rfloor}
\nu_j\in\mathbb N$ denote the digit-sum 
of a natural integer with binary expansion
$n=\sum_{j\geq 0}\nu_j2^j,\ \nu_0,\nu_1,\dots
\in\{0,1\}$.

\begin{thm} \label{thmmod2}
We have
$$\det(P(2n))=(-1)^n$$
and
$$\det(P(2n+1))=(-1)^{n+\mathop{ds}(n)}\ .$$
\end{thm}

\begin{rem}\label{remdetPF}
The infinite symmetric integral matrix $\tilde P$ with 
coefficients
$\tilde P_{s,t}={s+t\choose s}$ given by the binomial
coefficients is sometimes called the {\it Fermat matrix}.
Vandermonde's identity $\sum_{k=0}{s\choose k}{t\choose k}=
\sum_{k=0}{s\choose k}{t\choose t-k}={s+t\choose t}$ shows that
$\det(\tilde P(n))=1$ where $\tilde P(n)$ is the symmetric $n\times n$
submatrix with coefficients ${s+t\choose s},0\leq s,t<n$ of 
$\tilde P$.
\end{rem}

\subsection{$2-$valuations}\label{subsec2val}

Given a prime $p$, we denote by 
$v_p:\mathbb Q^*\longrightarrow\mathbb N$ the $p-$valuation.
Any rational number $\alpha$ can thus be written
in the form $\alpha=p^{v_p(\alpha)} \frac{n}{m}$ with 
$n,m\in\mathbb Z$ coprime to $p$.
Let $V(n)$ be the symmetric $n\times n$ matrix with coefficients 
$V_{s,t}\in\{\pm 1,\pm i\}$ given by
$$V_{s,t}=i^{v_2({s+t\choose s})},\ 0\leq s,t<n\ .$$

The $p-$valuation $v_p({s+t\choose s})$ of a binomial coefficient 
can be computed using a Theorem of Kummer
stating that $v_p({s+t\choose s})$ equals the number of carries 
occuring during the addition of the $p-$ary integers 
$s=\sum_{j\geq 0}\sigma_jp^j$ and $t=\sum_{j\geq 0}\tau_jp^j$.
More precisely, Kummer shows the identity
$$v_p(x!)=\frac{1}{p-1}\left(x-\sum_{j\geq 0}\xi_j\right)$$
(see Lehrsatz, page 115 of \cite{Kummer}) where $x=\sum_{j\geq 0}\xi_jp^j
\in\mathbb N$ with $\xi_j\in\{0,1,\dots,p-1\}$. This implies the
formula
$$v_p({s+t\choose s})=v_p((s+t)!)-v_p(s!)-v_p(t!)=
\frac{1}{p-1}\sum_{j\geq 0}(\sigma_j+\tau_j-u_j)$$
(see \cite{Kummer}, page 116) where $\sigma_j,\tau_j,u_j\in
\{0,1,\dots,p-1\}$ are defined by the $p-$ary expansions
$s=\sum_{j\geq 0}\sigma_jp^j,
t=\sum_{j\geq 0}\tau_jp^j$ and $s+t=\sum_{j\geq 0}u_jp^j$.

The next result uses the regular folding sequence
$f:\{1,2,\dots \}\longrightarrow\{\pm 1\}$. It is defined recursively by 
$f(2^n)=1$ and $f(2^n+a)=-f(2^n-a)$ for $1\leq a<2^n$, see 
for example \cite{AS}.

\begin{thm}\label{thmCRAS} We have
$$\det(V(2n))=(-1)^n\ \prod_{k=1}^{2n-1}(1-f(k)i)\in\mathbb Z[i]$$
and 
$$\det(V(2n+1))=(-1)^{n+\mathop{ds}(n)}\
\prod_{k=1}^{2n}(1-f(k)i)\in\mathbb Z[i]$$
(with $\mathop{ds}(\sum_{j=0}\nu_j2^j)
=\sum_{j\geq 0}\nu_j$ denoting the binary digit-sum).
\end{thm}

\begin{rem} Let $D$ denote the diagonal matrix with diagonal 
entries $i^{\mathop{ds}(0)},i^{\mathop{ds}(1)},i^{\mathop{ds}(2)},
\dots$. The paper \cite{BacherCRAS} deals with 
the Hankel matrix $H$ defined by $H_{s,t}=i^{ds(s+t)},
\ 0\leq s,t$ related by $H=D\overline V D$ to the complex
conjugate $\overline V$ of the matrix $V$ involved in Theorem
\ref{thmCRAS}.

Let us also mention that slight extensions of 
the computations occuring in our proof 
of Theorem \ref{thmCRAS} establish the existence of nice 
continued $J-$fraction expansions for the formal power series
(cf. \cite{BacherCRAS})
$$\begin{array}{l}
\displaystyle \prod_{k=0}^\infty(1+ix^{2^k})\ ,\\ 
\displaystyle \frac{1}{x}\left(\frac{1+i}{2}+\frac{1-x}{i-1}\prod_{k=0}^\infty(1+ix^{2^k})\right)\ ,\\
\displaystyle 
\frac{1}{x^2}\left(\frac{1+i}{2}+\frac{i-1}{2}x+\frac{1-x^2}{i-1}\prod_{k=0}^\infty(1+ix^{2^k})\right)\ .\end{array}$$
\end{rem}

\subsection{Beeblebrox reduction}

The idea (and the ``Beeblebrox'' terminology) of 
considering the ``Beeblebrox reduction'' of binomial coefficients are 
due to Granville, see \cite{Gr1} and \cite{Gr2}.

We define the ``{\it Beeblebrox reduction}''
as the Dirichlet character $\chi_B:\mathbb Z\longrightarrow\{0,\pm 1\}$ 
given by
$$\chi_B(x)=\left\lbrace\begin{array}{cl}
0&\hbox{if }x\equiv 0\pmod 2\\
1&\hbox{if }x\equiv 1\pmod 4\\
-1&\hbox{if } x\equiv 3\pmod 4\end{array}\right.$$
or equivalently by $\chi_B(2\mathbb Z)=0,\chi_B(4\mathbb Z\pm 1)=\pm 1$.
Beeblebrox reduction is the unique Dirichlet character modulo $4$
not factorising through $\mathbb Z/2\mathbb Z$.

The following result allows fast
computations of $\chi_B({n\choose k})$.

\begin{thm}\label{thmBeeblebroxformula} We have
$$\begin{array}{lcl}
\chi_B ({2n\choose 2k})&=&\chi_B ({n\choose k})\\
\chi_B ({2n\choose 2k+1})&=&0\\
\chi_B ({2n+1\choose 2k})&=&(-1)^k\ \chi_B ({n\choose k})\\
\chi_B ({2n+1\choose 2k+1})&=&(-1)^{n(k+1)}\ 
\chi_B ({n\choose k})\end{array}$$
\end{thm}

We denote by
$Z(n)$ (where the letter $Z$ stands for Zaphod Beeblebrox, following 
the amusing terminology of \cite{Gr1} and \cite{Gr2}) 
the symmetric {\it Beeblebrox matrix} of size $n\times n$
with coefficients $Z_{s,t}\in\{-1,0,1\}$ for $0\leq s,t<n$ given by
the Beeblebrox reduction $Z_{s,t}=\chi_B({s+t\choose s})$
of binomial coefficients.

Define $f:\mathbb N\longrightarrow \pm 3^{\mathbb Z}$
by $f(0)=1,f(1)=-1$ and recursively by
$$f(2^a+b)=\left\lbrace\begin{array}{ll} 3f(b)&\hbox{if }2b<2^a\\
\frac{1}{3}f(b)&\hbox{otherwise}\end{array}\right.$$
for $n=2^a+b\geq 2$ where $0\leq b<2^a$.
Laurent Bartholdi pointed out that the value $f(n)$
of a binary integer 
$n=\sum_{i\geq 0}\nu_i2^i$ is also given by
$$\begin{array}{lcl}
\displaystyle
f(n)&
\displaystyle=&
\displaystyle (-1)^n3^{\sharp\{i\ \vert\ 
\nu_i=0,\ \nu_{i+1}=1\}-
\sharp\{i\ \vert\ \nu_i=\nu_{i+1}=1\}}\\
&\displaystyle =&
\displaystyle (-1)^n
\prod_{i\geq 0}3^{(1-2\nu_i)\nu_{i+1}}\ .\end{array}$$

\begin{thm}\label{thmBeeblebroxdet} We have
$$\det(Z(n))=\prod_{k=0}^{n-1} f(k)\in \pm 3^{\mathbb N}\ .$$
\end{thm}

\subsection{The Jacobi-character modulo $8$}

Let $\chi_J:\mathbb N\longrightarrow\{0,\pm 1\}$ denote the Dirichlet
character modulo $8$ defined by $\chi_J(2\mathbb Z)=0$ and
$$\chi_J(n)=\left\lbrace\begin{array}{cl}
\displaystyle 1\qquad&\displaystyle\hbox{if }\ n\equiv \pm 1\pmod 8\ ,\\
\displaystyle -1\qquad&\displaystyle\hbox{if }\ n\equiv \pm 3\pmod 8\ .
\end{array}\right.$$
Since $\chi_J(p)\equiv 2^{(p-1)/2}\pmod 2$ for $p$ an odd prime, 
we call $\chi_J$ the Jacobi-character and 
we consider the matrix $J(n)$ with coefficients 
$J_{s,t}=\chi_J({s+t\choose s})$ for $0\leq s,t<n$. 
The techniques of this paper can be used to prove the following result:

\begin{thm}\label{thmJacDir} We have
$$\det(J(n))=\prod_{k=0}^{n-1}g(k)\in\pm 3^{\mathbb N}$$
where 
$$g(n)=(-1)^n\prod_{k\geq 0}3^{e(\lfloor n/2^k\rfloor)}$$
with $e(k)$ the $8-$periodic function given by the table
$$\begin{array}{|c|cccccccc|}
\hline
k&0&1&2&3&4&5&6&7\\
\hline
e(k)&0&0&1&-1&0&-2&3&-1\\
\hline\end{array}$$
\end{thm}

The matrices $Z(n)$ and $J(n)$ share many features.
However the considerable greater complexity of all objects
attached to $J$ makes it more difficult to pin down
interesting algebraic structures associated to $J$.

Let me also add that the remaining Dirichlet character modulo $8$
(given by $\tilde \chi(2\mathbb Z)=0$ and $\tilde \chi (\epsilon 
a+8\mathbb Z)=\epsilon$ for $a\in\{1,3\}$ and $\epsilon\in\{\pm 1\}$)
does not seem to give something interesting:
the symmetric infinite matrix obtained by 
applying $\tilde \chi$ to binomial coefficients has probably 
no $LU$ decomposition in the algebra of recurrence matrices.

Similarly, there seem to be no new interesting Dirichlet characters
$\pmod{16}$ (with values in the Gaussian integers $\{\pm 1,\pm \sqrt{-1}\}$.

\begin{rem} The case of the character $n\longmapsto \chi(n)=\epsilon_n\equiv 
n^{(p-1)/2}\pmod p$ with $\epsilon_n\in\{0,\pm 1\}$
for $p$ an odd prime gives rise to similar results which 
are somewhat trivial. Indeed, Lucas's theorem implies that 
the corresponding infinite matrix
with  coefficients $\chi({i+j\choose i}),\ 0\leq i,j$ 
has a structure of an infinite tensor-power (corresponding to a
``recurrence matrix of complexity $1$''). It is thus not very 
interesting and easy to handle. 

The same remark holds for 
the remaining characters $\pmod p$ when working in a suitable field
(or integral subring) of cyclotomic numbers.
\end{rem}

\subsection{$q-$binomials}

The expansion
$(x+y)^n=\sum_{k=0}^n{n\choose k}_qx^ky^{n-k}$
involving two non-commuting variables $x,y$ related 
by $yx=qxy$ where $q$ is a central variable defines the 
$q-$binomials coefficients
$${n\choose k}_q=\frac{\prod_{j=1}^n(1-q^j)}{\left(
\prod_{j=1}^k(1-q^j)\right)\left(
\prod_{j=1}^{n-k}(1-q^j)\right)}\in\mathbb N[q]\ .$$
An ordinary binary coefficient ${s+t\choose s}$ 
can be identified with the number of lattice paths with steps
$(1,0)$ and $(0,-1)$, starting at $(0,s)$ and ending at $(t,0)$.
Similarly, the coefficient of $q^c$ in the $q-$binomial ${s+t\choose s}_q$ 
counts the number of such paths delimiting a polygon of area 
$c$ in the first quadrant $\{(x,y)\in\mathbb R^2\ \vert\ x,y\geq 0\}$. 

Reflecting all paths contributing to ${s+t\choose s}_q$ with respect 
to the diagonal line $x=y$ yields the equality
$${s+t\choose s}_q={s+t\choose t}_q\ .$$

Rotating all paths contributing to ${s+t\choose s}_q$ by a half-turn
centered at $\frac{1}{2}(t,s)$ shows the identity
$${s+t\choose s}_q=q^{st}{s+t\choose s}_{q^{-1}}\ .$$

Partitioning all paths contributing to ${s+t\choose s}_q$
accordingly to the nature of their first step (horizontal or vertical)
shows the recursive formula
$${s+t\choose s}_q=q^s{s+t-1\choose s}_q+{s+t-1\choose s-1}_q$$
or equivalently ${n\choose k}_q=q^k{n-1\choose k}_q+{n-1\choose k-1}_q$
which is the $q-$version of the celebrated recurrence relation 
${n\choose k}={n-1\choose k}+{n-1\choose k-1}$ for 
ordinary binomial coefficients.

Cutting all lattice paths $\gamma$ contributing to ${s+t\choose s}_q$
along the diagonal line $s=t$ in two lattice paths shows the formula
$$\sum_k q^{k^2}{s\choose k}_q\ {t\choose k}_q={s+t\choose s}_q$$
where $k\in \{0,1,\dots,\mathop{min}(s,t)\}$.
This identity amounts to the matrix identity
$P_q=L_qD_qL_q^t$ where $P_q$ is the infinite symmetric matrix 
with coefficients ${s+t\choose s}_q,\ 0\leq s,t$, where
$L_q$ is the lower triangular
unipotent matrix with coefficients ${s\choose t}_q,\ 0\leq s,t$ and where
$D_q$ is diagonal
with diagonal coefficients $1,q,q^4,q^9,q^{16},q^{25},\dots$.
Denoting by $P_q(n)$ the submatrix ${s+t\choose s}_q,\ 0\leq s,t<n$
formed by the first $n$ rows and columns of $P_q$ we have the identity
$\det(P_q(n))=q^{\sum_{j=0}^{n-1}j^2}$ which specialises to 
the identity $\det(P_1(n))=1$ of Remark \ref{remdetPF}.
Appendix I of \cite{GR} contains many more formulae for ${n\choose k}_q$.

The formula of Lucas
$${a\choose b}\equiv {\lfloor a/p\rfloor\choose \lfloor
  b/p\rfloor}\ 
{a\pmod p\choose b\pmod p}\pmod p$$
(where $p$ is a prime number and
where $a\pmod p,b\pmod p\in\{0,1,\dots,p-1\}$),
see Section \ref{subsectmod2} or \cite{Lucas}, has the following known
analogue for $q-$binomials which reduces their evaluation at roots 
of $1$ of small order to evaluations of ordinary binomial coefficients.

\begin{thm} \label{thmformqbinn}
If $\omega=e^{2i\pi k/n}$ is a
primitive $n-$th root of $1$ (ie. $(k,n)=1$ with $k\in\mathbb Z$ and
$n\in\mathbb N$) then
$${a\choose b}_\omega={\lfloor a/n\rfloor\choose \lfloor
  b/n\rfloor}_1\ 
{a\pmod n\choose b\pmod n}_\omega$$
for all $a,b\in\mathbb N$ 
where $a\pmod n,b\pmod n\in\{0,1,\dots,n-1\}$. 
\end{thm}

Theorem \ref{thmformqbinn} can be used to establish formulae for 
determinants of 
the symmetric matrices obtained by considering the reduction modulo $2$,
the Beeblebrox reduction or the reduction using the Jacobi character 
modulo $8$ of (the real and imaginary part of)
${s+t\choose s}_q,\ 0\leq s,t<n$ evaluated at $q=-1$ and $q=i$.

\section{The algebra of recurrence matrices}\label{sectalgrec}

Recurrence matrices, introduced in \cite{BacherCRAS}, are a convenient 
tool for proving our main results.
Recurrence matrices are closely related 
to rational formal power series in free non-commutative variables
and can be considered as generalisations of finite state automata or 
of iterated tensor products. They arise also naturally in the context
of ``automata groups'', a notion generalising a famous group of Grigorchuk,
see \cite{Grigorchuk}. The following exposition does not strive for 
exhaustivity or for the largest possible generality. 
Generalisations (e.g. by replacing the field
of complex numbers by an arbitrary commutative field or
by considering 
sequences of square matrices of size $k^n\times k^n,n\in \mathbb N$ 
for $k\in\{1,2,3,\dots\}$) are fairly straightforward and
contained in \cite{BacherCRAS} or with more details in
\cite{BacherRM}.

The papers \cite{Bar} and \cite{Sid} deal with interesting subalgebras,
called self-similar algebras, formed by recurrence matrices.

\subsection{Recurrence matrices}

Consider the vector space 
$$\mathcal A=\prod_{n=0}^\infty M_{2^n\times 2^n}(\mathbb C)$$
whose elements are sequences $A=(A[0],A[1],A[2],\dots)$
with $A[n]\in M_{2^n\times 2^n}(\mathbb C)$ denoting a complex square
matrix of size $2^n\times 2^n$. The obvious product
$$AB=(A[0]B[0],A[1]B[1],A[2]B[2],\dots)$$
turns $\mathcal A$ into an associative algebra.
Denoting by 
$$\rho(0,0)A,\rho(0,1)A,\rho(1,0)A,\rho(1,1)A\in\mathcal A$$
the four ``corners'' of 
$$A=A[0],\left(\begin{array}{cc}(\rho(0,0)A)[0]&(\rho(0,1)A)[0]\\
(\rho(1,0)A)[0]&(\rho(1,1)A)[0]\end{array}\right),
\left(\begin{array}{cc}(\rho(0,0)A)[1]&(\rho(0,1)A)[1]\\
(\rho(1,0)A)[1]&(\rho(1,1)A)[1]\end{array}\right),\dots$$
obtained  (after deletion of the $1\times 1$ matrix $A[0]$)
by considering for all $n\geq 1$ 
the $2^{n-1}\times 2^{n-1}$ submatrix defined by
the first or last $2^{n-1}$ rows and by the first or last $2^{n-1}$
columns of $A[n]$,
we get four linear endomorphisms $\rho(s,t)\in\mathrm{End}(\mathcal A),
\ 0\leq s,t\leq 1$, of the vector space $\mathcal A$.
We call these endomorphisms {\it shift maps}.
Using a hopefully suggestive synthetic notation, an element 
$A\in\mathcal A$ can thus be written as
$$A=A[0],\left(\begin{array}{cc}\rho(0,0)A&\rho(0,1)A\\\rho(1,0)A&\rho(1,1)A
\end{array}\right)$$
with $A[0]\in \mathbb C$ and  
$\rho(0,0)A,\rho(0,1)A,\rho(1,0)A,\rho(1,1)A\in \mathcal A$.

{\bf Definition} A subspace $\mathcal V\subset \mathcal A$ is {\it
recursively closed} if $\rho(s,t)\mathcal V\subset \mathcal V$
for all $s,t$.

The {\it recursive closure} $\overline{\mathcal S}$ of a subset 
$\mathcal S\in\mathcal A$ is 
the smallest recursively closed subspace of $\mathcal A$ which contains
$\mathcal S$. We denote by $\overline A$ the recursive closure of
the subset $\{A\}$ reduced to a single element $A\in\mathcal A$.
The {\it complexity} of $A\in\mathcal A$ is the dimension 
$\mathop{dim}(\overline A)\in
\mathbb N\cup\{\infty\}$ of the recursive closure
$\overline A\subset \mathcal A$.

An element $A\in\mathcal A$ is a {\it recurrence matrix} 
if its recursive closure 
$\overline A$ is of finite dimension. 
We denote by $\mathcal R\subset \mathcal A$ the subset of
all recurrence matrices.

Writing $\rho(X_{s,t})A$ or simply $X_{s,t}A$ 
for $\rho(s,t)A,\ 0\leq s,t\leq 1$, 
the shift maps $\rho(s,t)\in\mathrm{End}(\mathcal A)$ 
induce a linear representation (still denoted)
$\rho:\{X_{0,0},X_{0,1},X_{1,0},X_{1,1}\}^*\longrightarrow
\mathrm{End}(\mathcal A)$, recursively defined by
$$(X_{s_1,t_1}X_{s_2,t_2}\cdots X_{s_l,t_l})A=
(X_{s_1,t_1}X_{s_2,t_2}\cdots X_{s_{l-1},t_{l-1}})(\rho(s_l,t_l)A)\ ,$$
of the free non-commutative monoid 
$\{X_{0,0},X_{0,1},X_{1,0},X_{1,1}\}^*$, called the
{\it shift monoid}, 
in four generators $X_{0,0},X_{0,1},X_{1,0},X_{1,1}$ representing 
shift maps. Subrepresentations of $\rho$ correspond to 
recursively closed subspaces 
$\mathcal V$ of $\mathcal A$ spanned by (unions of) orbits 
under $\{X_{0,0},X_{0,1},X_{1,0},X_{1,1}\}^*$. 

The linear action of the monoid $\{X_{0,0},X_{0,1},X_{1,0},X_{1,1}\}^*$ 
on $\mathcal A$ suggests to consider the bijective map 
which associates an element $A\in\mathcal A$ with the 
non-commutative formal power series
$$\sum_{\mathbf X\in\{X_{0,0},X_{0,1},X_{1,0},X_{1,1}\}^*}
((\mathbf X A)[0])\mathbf X\in\mathbb C\langle\!\langle
X_{0,0},X_{0,1},X_{1,0},X_{1,1}\rangle\!\rangle$$
in four free non-commutative variables $X_{0,0},X_{0,1},X_{1,0},X_{1,1}$.
This bijection restricts to a bijection between the vector
space $\mathcal R$ of recurrence matrices and rational
elements in $\mathbb C\langle\!\langle
X_{0,0},X_{0,1},X_{1,0},X_{1,1}\rangle\!\rangle$.

The algebraic structure of $\mathcal R\subset \mathcal A$
is described by the following result.

\begin{prop} \label{propadmult}
(i) We have $\mathop{dim}(\overline{\lambda A})=
\mathop{dim}(\overline A)$ for all $\lambda\in\mathbb C^*$ and for all
$A\in\mathcal A$.

\ \ (ii) We have $\mathop{dim}(\overline{A+B})\leq
\mathop{dim}(\overline A)+\mathop{dim}(\overline B)$
for all $A,B\in\mathcal A$.

\ \ (iii) We have  $\mathop{dim}(\overline{AB})\leq
\mathop{dim}(\overline A)\mathop{dim}(\overline B)$
for all $A,B\in\mathcal A$.
\end{prop}

\begin{rem}
The inequalities of assertion (ii) and (iii) can of course be strict:
Consider two elements $A,B\in \mathcal A$ defined by $A[n]=\frac{1+(-1)^n}{2}
\mathrm{Id}[n],\ B[n]=\frac{1-(-1)^n}{2}\mathrm{Id}[n]$
where $\mathrm{Id}[n]$ denotes the identity
matrix of size $2^n\times 2^n$. The elements $A,$B have
common recursive closure
$\overline A=\overline B=\mathbb C A+\mathbb C B$  of dimension $2$.
Their sum $A+B=\mathrm{Id}\in\mathcal R$ is the identity element 
having complexity $1$ and their product $AB=0$ has complexity $0$.
\end{rem}

\begin{cor} The set $\mathcal R$ of recurrence matrices
is a subalgebra of $\mathcal A$.
\end{cor} 

\noindent
{\bf Proof of Proposition \ref{propadmult}} (i) and (ii) are obvious. 

Denoting (slightly abusively) by 
$\overline A\ \overline B=\{\sum X_iY_i\ \vert X_i\in\overline 
A,Y_i\in\overline B\}$
the vector space spanned by all products $XY,\ X\in\overline
A,Y\in\overline B$, we have
$AB\in \overline A\ \overline B$.

For
$$XY=(X[0]Y[0]),\left(\begin{array}{cc}
\rho(0,0)(XY)&\rho(0,1)(XY)\\
\rho(1,0)(XY)&\rho(1,1)(XY)\\
\end{array}\right)
\in \overline A\ \overline B$$
with $X\in\overline A,Y\in\overline B$, the computation
$$\begin{array}{l}
\rho(0,0)(XY)=(\rho(0,0)X)(\rho(0,0)Y)+(\rho(0,1)X)(\rho(1,0)Y)\\
\rho(0,1)(XY)=(\rho(0,0)X)(\rho(0,1)Y)+(\rho(0,1)X)(\rho(1,1)Y)\\
\rho(1,0)(XY)=(\rho(1,0)X)(\rho(0,0)Y)+(\rho(1,1)X)(\rho(1,0)Y)\\
\rho(1,1)(XY)=(\rho(1,0)X)(\rho(0,1)Y)+(\rho(1,1)X)(\rho(1,1)Y)
\end{array}$$
shows that $\overline A\ \overline B$ is recursively closed
of dimension $\leq 
\mathop{dim}(\overline A)\mathop{dim}(\overline B)$.
Assertion (iii) follows now from the obvious inclusion 
$\overline{AB}\subset \overline A\ \overline B$.
\hfill$\Box$

\begin{rem} Certain properties of binomial coefficients $\pmod{p^d}$
are easy to study using ``recurrence matrices'' 
given by sequences of matrices of size $p^j\times p^j,j=0,1,\dots$
with entries in the associative ring $\mathbb Z/
p^d\mathbb Z$.
\end{rem}

\subsection{Recursive presentations}

An element $A\in\mathcal A$ is completely determined
by the action of the shift maps $\rho(s,t)$ on its recursive 
closure $\overline A$ (spanned
by $\{X_{0,0},X_{0,1},X_{1,0},X_{1,1}\}^*A$),
together with the restriction to $\overline A$ of the {\it augmentation
map} $\pi_0\in\mathcal A^*$ defined by projecting an
element $X=(X[0],X[1],\dots)\in\mathcal A$ onto its {\it initial value}
$\pi_0(X[0],X[1],X[2],\dots)=X[0]\in\mathbb C$.

A recurrence matrix $A$ can thus be given by
a finite amount of data: An (expression describing the) 
element $A$ of the finite-dimensional
vector space $\overline A$, the restriction 
(still denoted) $\pi_0\in \overline A^*$ 
of the augmentation map expressing the initial values of elements
in $\overline A$, and a $2\times 2$ matrix $\rho=\left(\begin{array}{cc}
\rho(0,0)&\rho(0,1)\\ \rho(1,0)&\rho(1,1)\end{array}\right)
\in M_{2\times 2}(\overline A\otimes \overline A^*)$ of tensors
encoding the shift maps. The coefficients of the matrix 
$A[n]\in M_{2^n\times 2^n}$ are then obtained by ``contractions'' of 
$\pi_0\ \rho^n\  A$.

This leads to the notion of {\it recursive presentations}. 
A recursive presentation for $A\in\mathcal R$ is given 
by the choice of a basis $A_1=A,\dots,A_a$
of $\overline A$ (or more generally of a finite set $A_1,\dots,A_a$
spanning a recursively closed 
vector space containing $\overline A$) and by recursive 
identities 
$$A_j=\pi_0(A_j),\left(\begin{array}{cc}
\rho(0,0)A_j=\sum_{k=1}^a\rho(0,0)_{k,j}A_k&
\rho(0,1)A_j=\sum_{k=1}^a\rho(0,1)_{k,j}A_k\\
\rho(1,0)A_j=\sum_{k=1}^a\rho(1,0)_{k,j}A_k&
\rho(1,1)A_j=\sum_{k=1}^a\rho(1,1)_{k,j}A_k\end{array}\right)$$
encoding the initial values $\pi_0(A_j)$ and the values 
$\rho(s,t)(A_j)\in\overline A$ of the shift maps. A recursive 
presentation defines the elements $A_1=A,\dots,A_a$ spanning
(a recursively closed subspace containing) 
$\overline A$ recursively by expressing the 
four ``blocks'' of $A_j[n+1]$ as linear combinations of 
$A_1[n],\dots,A_a[n]$.

\begin{rem} We use the convention that $0\in\mathcal R$ admits
the empty presentation and $0$ is ``the'' element of an
empty basis. We speak thus of ``the'' basis $A_1,\dots$
of $\overline 0$ representing $0=A_1$.
\end{rem}

\subsection{Saturation level}\label{sectsatlev}

We denote by $\pi_l(A)=A[l]$ the projection of a matrix sequence
$A=(A[0],A[1],A[2],\dots)\in\mathcal A$ onto its
square matrix $A[l]$ of size $2^l\times 2^l$. Similarly, 
$$\pi_{\leq l}(A)=(\pi_0(A),\pi_1(A),\dots,\pi_l(A))=
(A[0],A[1],\dots,A[l])\in\oplus_{j=0}^l M_{2^j\times 2^j}(\mathbb C)$$
denotes the projection of the sequence $A$ onto its first $l+1$
matrices. 

The {\it saturation level} of a finite dimensional subspace
$\mathcal V\subset \mathcal A$ is the smallest integer $N\in\mathbb N$
such that $\mathcal K_{\leq N}(\mathcal V)=\mathcal K_{\leq N+1}
(\mathcal V)$ where 
$\mathcal K_{\leq l}(\mathcal V)\subset \mathcal V$ is the kernel 
of the projection 
$\pi_{\leq l}:\mathcal V\longrightarrow
\pi_{\leq l}(\mathcal A)=\oplus_{j=0}^l M_{2^j\times 2^j}$.

\begin{prop} \label{propsatlevel} We have 
$\mathcal K_{\leq N}(\mathcal V)=\{0\}$ for the saturation level $N$ of a
finite-dimensional subspace $\mathcal V\subset A$ which is recursively
closed. 

In particular, $\pi_{\leq N}:\mathcal V\longrightarrow
\oplus_{j=0}^N M_{2^j\times 2^j}$ defines an injection.
\end{prop}

\noindent
{\bf Proof} The obvious inclusions 
$\rho(s,t)\mathcal K_{\leq l+1}(\mathcal V)\subset \mathcal K_{\leq l}
(\mathcal V)$
imply that $\mathcal K_{\leq N}(\mathcal V)=\mathcal K_{\leq N+1}(\mathcal V)
\subset\mathcal V$
is recursively closed. Since the restriction to 
$\mathcal K_{\leq N+1}(\mathcal V)\subset \mathcal K_{\leq 0}$
of the augmentation map $\pi_0:\mathcal A\longrightarrow\mathbb C$
is trivial, we have $(\mathbf X K)[0]=0$ for all 
$\mathbf X\in\{X_{0,0},X_{0,1},X_{1,0},X_{1,1}\}^*$ and for all
$K\in\mathcal K_{\leq N}(\mathcal V)$. This shows
$\mathcal K_{\leq N}(\mathcal V)=\{0\}$.\hfill$\Box$

Proposition \ref{propsatlevel} 
enables us to extract a basis from a finite set 
$\mathcal S$ spanning a recursively closed vector space $
\mathcal V\subset \mathcal R$. Similarly, Proposition 
\ref{propsatlevel} allows the construction of a basis 
of the subspace $\overline A\subset \mathcal V$ for an element 
$A\in \mathcal V$ of a finite-dimensional recursively closed 
vector space $\mathcal V\subset \mathcal R$.

These operations are the necessary ingredients for effective
computations in the algebra $\mathcal R$. Effectivity means
that there exists an algorithm involving only a finite
number of elementary operations in the groundfield $\mathbb C$
and a finite amount of data which computes the result
of an algebraic expression (given by a non-commutative 
polynomial) involving (recursive 
presentations of) a finite number of elements in $\mathcal R$.

The necessary elementary algorithms can be briefly described 
as follows:

\subsubsection{Multiplication of $A\in\mathcal R$ by a non-zero scalar
$\lambda\in\mathbb C^*$}\label{subsubscal} A presentation of 
$\lambda A$ is obtained 
from a presentation of $A$
by multiplying the initial values $\pi(A_j)=A_j[0]\in\mathbb C$ 
with $\lambda$ (and by keeping the same shift maps).

\subsubsection{Addition of two elements $A,B\in\mathcal R$}
\label{subsubadd}
Add a first element $A_1+A_2$ 
having the obvious initial value $\pi_0(A+B)=A_1[0]+B_1[0]$
to the list of not necessarily linearly independent
elements $A_1,\dots,A_a,B_1,\dots,B_b$ 
spanning $\overline{A+B}$. The elements $\rho(s,t)(A_1+A_2),\rho(s,t)A_j,
\rho(s,t)B_j$ 
are given by 
$$\rho(s,t)(A_1+B_1)=\sum_{k=1}^a \rho(s,t)^A_{k,1}A_k+\sum_{k=1}^b
\rho(s,t)^B_{k,1}B_k$$
and 
$$\rho(s,t)A_i=\sum_{k=1}^a\rho(s,t)^A_{k,i}A_k,\ 
\rho(s,t)B_j=\sum_{k=1}^b\rho(s,t)^B_{k,j}B_k$$
for all $0\leq s,t\leq 1,\ 1\leq i\leq a,\ 1\leq j\leq b$
where $\rho(s,t)^A$ and $\rho(s,t)^B$ are the obvious shift
maps with respect to bases $A=A_1,\dots,A_a$ 
and $B=B_1,\dots,B_b$ of $\overline A$ and $\overline B$.
Working with the finite-dimensional recursively
closed subspace  $\mathbb C(A_1+B_1)+\sum_{i=1}^a
\mathbb C A_i+\sum_{j=1}^b\mathbb C B_j\subset\mathcal R$ 
one can now give a presentation of $A_1+B_1$ by
computing a basis of $\overline{A_1+B_1}$, followed by the computation
of the coefficients (with respect to this basis) of the shift maps.

\begin{rem}
The algorithms \ref{subsubscal} and \ref{subsubadd} can be used
to compare two elements $A,B\in\mathcal R$ by
computing a presentation of $A-B$.
\end{rem}

\subsubsection{Multiplication of two elements $A,B\in\mathcal R$}
\label{subsubmult}

Consider the $ab$ elements $C_{i,j}=A_iB_j,\ 1\leq i\leq a,\ 
1\leq j\leq b$ with initial values $C_{i,j}[0]=\pi_0(A_iB_j)=
A_i[0]B_j[0]$. Shift maps are given by
$$\leqno{(*)}
\rho(s,t)C_{i,j}=\sum_{k=1}^a\sum_{l=1}^b
\left(\rho(s,0)^A_{k,i}\rho(0,t)^B_{l,j}+
\rho(s,1)^A_{k,i}\rho(1,t)^B_{l,j}\right)C_{k,l}$$
using the notations of \ref{subsubadd}.
One constructs now a recursive presentation of $C_{1,1}=AB$
by proceeding as above using the recursively closed vector space 
$\sum_{i=1}^a\sum_{j=1}^b 
\mathbb CC_{i,j}\supset \overline{C_{1,1}}$.

\begin{rem}\label{rembirec}
The formulae (*) occuring in \ref{subsubmult}
define an associative product on the set $\mathcal L$ of all 
equivalence-classes of finite-dimensional linear representations 
of the monoid $\{X_{0,0},X_{0,1},X_{1,0},X_{1,1}\}^*$.
Considering also direct sums of
linear representations turns $\mathcal L$ into a semi-ring.
The semi-ring $\mathcal L$ has a homomorphism $\varphi$
into the semi-ring 
(with addition $\mathcal V+\mathcal W=\{X+Y\ \vert\ X\in\mathcal V,
Y\in\mathcal W\}$ and 
product $\mathcal V\mathcal W=\{\sum_i X_iY_i\ \vert\ X_i\in\mathcal V,
Y_i\in\mathcal W\}$)
formed by all finite-dimensional recursively closed subspaces
of $\mathcal R$. The elements of $\varphi(L)$ 
have two equivalent descriptions:
They can be identified with the subset $\mathcal L'
\subset \mathcal L$ given by all equivalence classes of finite-dimensional
linear representations of $\{X_{0,0},X_{0,1},X_{1,0},X_{1,1}\}^*$
containing no equivalence class of 
a subrepresentation with multiplicity $>1$.
The second descriptions involves birecursively closed vector spaces
which are defined as follows:
A vector space $\mathcal V\subset \mathcal R$ is 
{\it birecursively closed} if
it is recursively closed and if an arbitrary generic modification 
of the initial values in a recursive presentation of an element 
in $\mathcal V$ yields again a presentation of an element in $\mathcal V$.

The homomorphism of semi-rings $\varphi$ associates
to a finite-dimensional linear representation $\rho_f$ of
$\{X_{0,0},X_{0,1},X_{1,0},X_{1,1}\}$ the unique maximal birecursively
closed subspace of $\mathcal R$ whoses shift maps involve only
equivalence classes of subrepresentations in $\rho_f$.

Finite-dimensional birecursively closed subspaces of $\mathcal R$
are stable under direct sums and products and their semi-ring 
is the quotient semi-ring $\varphi(\mathcal L)$ of $\mathcal L$.
\end{rem}

\subsection{The $LU$ decomposition of a convergent non-singular
element in $\mathcal R$}\label{subsecLU}

An element $P\in\mathcal A$ such that $P=\rho(0,0)P$ is called 
{\it convergent}. It is given by considering 
the sequence 
$$P_{0,0},\left(\begin{array}{cc}P_{0,0}&P_{0,1}\\P_{1,0}&P_{1,1}\end{array}\right),\left(\begin{array}{cccc}P_{0,0}&P_{0,1}&P_{0,2}&P_{0,3}\\
P_{1,0}&P_{1,1}&P_{1,2}&P_{1,3}\\
P_{2,0}&P_{2,1}&P_{2,2}&P_{2,3}\\
P_{3,0}&P_{3,1}&P_{3,2}&P_{3,3}\end{array}\right),\dots$$
of all square submatrices formed by the first $2^n$ rows and columns
of an infinite ``limit''matrix 
$$\left(\begin{array}{ccccc}P_{0,0}&P_{0,1}&P_{0,2}&\dots\\
P_{1,0}&P_{1,1}&P_{1,2}&\dots\\
\vdots\end{array}\right)\ .$$
Henceforth we denote generally a convergent element in $\mathcal A$ and the
associated infinite matrix by the same letter. This should
not lead to confusions except in cases where both interpretations 
are correct.

We call an infinite matrix $P$ {\it non-singular}
if the $k\times k$ square matrix $P(k)$ formed by its first $k$ rows and
columns has non-zero determinant for all $k\geq 1$. Such a non-singular
matrix $P$ has an $LU-$decomposition: It can be written as 
$P=LU$ with $L$ lower triangular unipotent ($1$'s on the diagonal)
and $U$ upper triangular non-singular. The identity $P=LU$ 
implies the equality $\det(P(k))=\det(U(k))$ for all $k\geq 1$
and gives rise to an $LU-$decomposition 
in $\mathcal A$ by considering as above for $n=0,1,2,\dots$
the submatrices formed by
the first $2^n$ rows and columns of of $P,L$ and $U$. If
$P$ is symmetric we have moreover $U=DL^t$ where $D$ is diagonal
non-singular and $L^t$ is obtained by transposing the matrix $L$.

All proofs of the results presented in Section \ref{sectMain} boil
down to $LU-$decompositions with $P=P^t=LU,L,D,U=DL^t\in\mathcal R$.

\begin{rem} Call an element $A\in\mathcal A$ non-singular
if it involves only non-singular matrices $A[0],A[1],\dots$. 
Such an element
has an $LU-$decomposition (in the obvious sense) in $\mathcal A$.

Proving the non-existence of an $LU-$decomposition in $\mathcal R$
for a suitable given non-singular recurrence matrix $A\in\mathcal R$ 
is probably difficult. 

The related problem of constructing the (existing) 
recurrence matrices $L,U\in\mathcal R$ 
from the knowledge (of a recursive presentation) of $A=LU\in\mathcal R$
has however an algorithmic answer: One
proceeds as for the existence of an inverse element by guessing
recursive presentations for $L$ and $U$ using $LU-$decompositions
of finitely many matrices $A[0],A[1],\dots,A[N+1]$.
In case of succes, the resulting hypothetical decomposition, 
if correct, can then be proven 
to hold. The necessary algorithm is obtained after minor 
modifications from the algorithm for 
computing the inverse $A^{-1}\in\mathcal R$ of an invertible
element $A\in\mathcal R$ described in Section \ref{sectalgoinv}
\end{rem}

\section{Invertible recurrence matrices}\label{sectinv}

The set of all recurrence matrices which are invertible in
the algebra $\mathcal R$ forms the group of {\it units}
in $\mathcal R$. Determining the inclusion in the unit group
of $\mathcal R$ of 
a recurrence matrix $A$ is perhaps a difficult problem without 
algorithmic solution. Indeed, we have the following result.

\begin{prop} \label{propinv} 
(i) For every natural integer $n$ there exist invertible
recurrence matrices $A,B=A^{-1}\in\mathcal R$ such that $\dim(
\overline{A})=2$ and $\dim(\overline B)> n$. 

\ \ (ii) There exist elements in $\mathcal R$ which 
are invertible in the algebra $\mathcal A$ but not in the 
subalgebra $\mathcal R$ of recurrence matrices.
\end{prop}

\begin{rem} The assumption $\dim(\overline A)=2$ is optimal:
Invertible recurrence matrices of complexity
$1$ form a subgroup (isomorphic to $\mathbb C^*\times 
\hbox{GL}_2(\mathbb C)$) in $\mathcal R$.
\end{rem}

\noindent
{\bf Proof of Proposition \ref{propinv}} 
For $\omega\in\mathbb C^*$, consider the convergent element
$$A=1,\left(\begin{array}{cc}1\\
-\omega&1\end{array}\right),\left(\begin{array}{cccc}1\\
-\omega&1\\0&-\omega&1\\0&0&-\omega&1\end{array}\right),\dots\in\mathcal A$$ 
consisting of lower triangular unipotent matrices with constant subdiagonal
$-\omega$. It defines a recurrence matrix $A=A_1$ of complexity $2$ 
recursively presented by
$$
A_1=1,\left(\begin{array}{cc}A_1&0\\ A_2&A_1\end{array}\right),
\quad 
A_2=-\omega,\left(\begin{array}{cc}0&A_2\\0&0\end{array}\right)\ .$$
%Setting $A_1=A$ and 
%$$A_2=\rho(1,0)A_1=-\omega,\left(\begin{array}{cc}0&-\omega\\
%0&0\end{array}\right),\left(\begin{array}{cccc}
%0&0&0&-\omega\\0&0&0&0\\0&0&0&0\\0&0&0&0\end{array}\right),\dots$$
%we have
%$$\begin{array}{l}
%\displaystyle
%\rho(0,0)A_1=\rho(1,1)A_1=A_1,\ \rho(0,1)A_1=0,\ \rho(1,0)A_1=A_2\ ,\\
%\displaystyle
%\rho(0,0)A_2=\rho(1,0)A_2=\rho(1,1)A_2=0,\ \rho(0,1)A_2=A_2\ .
%\end{array}$$
%This shows that $A$ is a recurrence matrix of complexity $2$.

Since $A$ is given by a sequence of unipotent lower triangular matrices,
it is invertible in the algebra $\mathcal A$ with inverse
the convergent element
$$B=A^{-1}=1,\left(\begin{array}{cc}1\\
\omega&1\end{array}\right),\left(\begin{array}{cccc}1\\
\omega&1\\ \omega^2&\omega&1\\ \omega^3&\omega^2&\omega&1
\end{array}\right),\dots\in\mathcal A$$
whose limit is the infinite unipotent lower triangular
Toeplitz matrix with constant subdiagonals
associated to the geometric progression $1,\omega,\omega^2,\dots$.

For $k$ a strictly positive integer, we consider the element
$$S_k=\omega^{2^0k},\omega^{2^1k}\left(\begin{array}{cc}
1&\omega\\ \omega&\omega^2\end{array}\right),
\omega^{2^2k}\left(\begin{array}{cccc}
1&\omega&\omega^2&\omega^3\\\omega&\omega^2&\omega^3&\omega^4\\
\omega^2&\omega^3&\omega^4&\omega^5\\
\omega^3&\omega^4&\omega^5&\omega^6\end{array}\right),\dots\in\mathcal A$$
given by $\left(S_k[n]\right)_{i,j}=\omega^{2^nk+i+j},\ 0\leq i,j<2^n$.
A straightforward computation shows $\rho(s,t)S_k=S_{s+t+2k}$
for all $s,t$. Since we have 
$$B=1,\left(\begin{array}{cc}B&0\\S_1&B\end{array}\right)$$
(using the notation of recursive presentations),
%$$\rho(0,0)B=\rho(1,1)B=B,\ \rho(0,1)B=0,\ \rho(1,0)B=S_1\ ,$$
the algebraic closure $\overline B$ of $B$ is spanned by the set
$$\{X_{0,0},X_{0,1},X_{1,0},X_{1,1}\}^*B=B\cup\{S_k\ \vert\ 
k\geq 1\}\ .$$
Since $\omega^N=1$ implies $S_{k+N}=S_k$,
this set contains at most $1+N$ elements 
if $\omega$ is a root of $1$ of finite order $N$.

In order to prove assertion (i), we consider
the case where $\omega$ is a root of $1$ having odd order $N>2^n$ and 
we denote by $a\geq \mathop{log}_2(N+1)>n$ the order of $2$ in the 
multiplicative group $(\mathbb Z/N\mathbb Z)^*$.
The action of the monoid $\rho(0,0)^{\mathbb N}$ on $S_k$ 
corresponds then to the action of the Galois map $\omega\longmapsto 
\omega^2$ on the upper left coefficients 
$$\omega^{2^0k},\omega^{2^1k},\omega^{2^2k},\omega^{2^3k},\dots$$
of $S_k$. Elementary number theory (for instance reduction modulo
$2$ by choosing $\omega$ among the primitive $N-$th roots of 
$1$ in the field extension $\mathbb F_{2^a}$ of degree $a$ over 
$\mathbb F_2$)
shows that the set $S_{2^{\mathbb N}}$ spans a subspace 
of dimension $a>n$ in $\overline B$. This implies assertion (i).

%Assertion (i) follows from the inequality
%$a\geq \frac{\mathop{log}(N+1)}{\mathop{log}\ 2}$
%for the order $a$ of the cyclic subgroup generated by $2$ in
%$\mathbb Z/N\mathbb Z$.
%This proves assertion (i).

Consider now $A$ as above with $\omega\in \mathbb 
C\setminus\overline{\mathbb Q}$ transcendental. 
This allows to consider $\omega$ as a variable and we have 
$\dim(\overline{B_\omega})\geq \dim(\overline{B_\xi})$ 
for the complexity of the inverse element 
$B_\omega=A_\omega^{-1}$ where $\xi\in\overline{\mathbb Q}$ is any algebraic
specialisation of $\omega$.
Assertion (ii) follows now from assertion (i).
\hfill $\Box$

\begin{rem} An example of $B=A^{-1}\in\mathcal A\setminus\mathcal R$
with $A\in\mathcal R$ invertible only in $\mathcal A$
is also given by $A,\ B$ as above with $\omega=2$
(the argument below works in fact for any $\omega$ of norm 
$\vert \omega\vert\geq 1$).
Indeed, otherwise, up to a constant, 
the initial values $2,2^2,2^4,2^8,\dots$
of the sequence $\rho(1,0)B,\rho(0,0)\rho(1,0)B,\rho(0,0)^2\rho(1,0)B,\dots$
should be bounded above by a geometric progression $1,\mu,\mu^2,\dots$
for any positive $\mu$ exceeding the spectral radius (absolute
value of the largest eigenvalue) of
$\rho(0,0)\in\mathrm{End}(\overline B)$.

The simplest element of $\mathcal R$ with an inverse in $\mathcal A\setminus
\mathcal R$ is perhaps given by the central diagonal 
recurrence matrix $A$ of complexity $2$
defined by $A[n]=(n+1)\hbox{Id}[n]$ where $\hbox{Id}[n]$ denotes the identity
matrix of size $2^n\times 2^n$. We give two proofs that
the  inverse element $A^{-1}$, given by $A^{-1}[n]=\frac{1}{n+1}
\hbox{Id}[n]$, has infinite complexity.

A first proof follows from the observation that the  
sequences of diagonal coefficients $\frac{1}{n+k+1}=(\rho(0,0)^kA^{-1})[n]$
form the rows of the Hilbert matrix of infinite rank with
coefficients $H_{i,j}=\frac{1}{1+i+j},\ 0\leq i,j$.

In order to give a second proof, we observe that all coefficients
of $A^{-1}$ are rational numbers and every prime number appears
as a denominator in a suitable coefficient of $A^{-1}$. This is impossible
for an element $B\in\mathcal R$ with rational coefficients. Indeed,
such an element $B$ has a recursive presentation with data 
(consisting of initial values and coefficients of shift maps with 
respect to a basis $\subset \{X_{0,0},X_{0,1},X_{1,0},X_{1,1}\}^*
B$ of $\overline B$) given by a finite set $\mathcal D$ of rational numbers. 
Since all coefficients of $B$ are evaluations of integral polynomials
on $\mathcal D$, all denominators of coefficients in $B$ involve 
only prime numbers occuring in the denominators of the finite set 
$\mathcal D$.
\end{rem}

\begin{rem} The quotient group of invertible elements modulo 
$\mathbb C^* \hbox{Id}$ can be turned into a metric group
by considering the positive real function 
$$A\longmapsto \parallel A\parallel=\max (\log(\dim(\overline A+\mathbb C\hbox{Id})),
\log(\dim(\overline{A^{-1}}+\mathbb C\hbox{Id})))$$
on the group $\Gamma$ of units in $\mathcal R$.
It satisfies $\parallel A\parallel \geq 0$ with $\parallel A\parallel=0$ 
only for $A\in\mathbb C^*\hbox{Id}$,
$\parallel AB\parallel\leq \parallel A\parallel +\parallel B\parallel$
and $\parallel A\parallel =\parallel A^{-1}\parallel$
and defines thus a left-invariant distance on the quotient group 
$\Gamma/\mathbb C^* \hbox{Id}$ by considering
$d(A,B)=\parallel A^{-1} B\parallel$.

The corresponding group over a finite field has finitely many elements
in balls of finite radii and it would be interesting to understand
the generating function
$$\sum_{A\in\Gamma(\mathbb F_{p^e})}t^{
\max(\dim(\overline A+\mathbb C\hbox{Id}),
\dim(\overline{A^{-1}}+\mathbb C\hbox{Id}))}\in\mathbb N[[t]]$$
where the sum is over all elements of the unit group $\Gamma(\mathbb 
F_{p^e})$ of the algebra $\mathcal R(\mathbb F_{p^e})$ 
defined in the obvious way over the finite field $\mathbb F_{p^e}$.

The related generating function 
$$\sum_{A\in\mathcal R(\mathbb F_{p^e})}t^{\dim(\overline A)}
\in\mathbb N[[t]]$$
counting all elements of given complexity in the algebra $\mathcal R(
\mathbb F_{p^e})$ is probably fairly  easy to compute.
It has convergency radius $0$ and is thus transcendental.
\end{rem}

\subsection{Algorithm for computing $A^{-1}\in\mathcal R$ for
$A\in\mathcal R$ invertible in $\mathcal R$}\label{sectalgoinv}

The following algorithm computes the inverse of an element 
$A\in\mathcal R$ if it has an inverse in $\mathcal R$ 
and fails (does never stop and uses more and more memory)
for an element $A$ as in assertion (ii) of 
Proposition \ref{propinv}. Non-invertibility in $\mathcal A$
will eventually be detected (assuming exact arithmetics over the ground
field) by exhibiting an integer $k$ for
which $A[k]$ is singular.

\noindent{\bf Input:} A (presentation of a) recurrence matrix $A\in\mathcal R$.

Set $N=1$ and $M=1$.

\noindent{\bf Loop} Compute, if possible, the matrices 
$B[0]=A[0]^{-1},B[1]=A[1]^{-1},\dots,
B[N+M+1]=A[N+M+1]^{-1}$.

If a matrix $A[k]$ with $k\leq N+M+1$ is not invertible,
print ``The matrix $A[k]$ is singular and $A$ has thus 
no inverse in $\mathcal A$" and stop.

We denote by $\tilde B\in\mathcal A$ the sequence $B[0],\dots,B[N+M+1]$
completed by arbitrary matrices $\tilde B[k]$ of size $2^k\times 2^k$ 
(for example by zero matrices) if $k>N+M+2$.

For $n,m\in\mathbb N$ such that $n+m\leq N+M+1$, we denote by 
$\mathcal V(n,m)\subset \oplus_{j=0}^{n} M_{2^j\times 2^j}(\mathbb C)$ 
the vector-space spanned by all $1+4+\dots+4^m$ elements
$(\pi_{\leq n}(\mathbf X \tilde B))_{\mathbf X\in\mathcal X^{\leq m}}$ 
where $\mathcal X^{\leq m}$ denotes the set of all $\frac{1-4^{m+1}}{1-4}=
1+4+4^2+\dots+4^m$ words of length $\leq m$ in the alphabet $\mathcal X=
\{X_{0,0},X_{0,1},X_{1,0},X_{1,1}\}$.

If $\dim(\mathcal V(N,M))<\dim(\mathcal V(N+1,M))$,
then increase $N$ by $1$ and iterate the Loop.

If $\dim(\mathcal V(N,M))<\dim(\mathcal V(N,M+1))$,
then increase $M$ by $1$ and iterate the Loop.

Setting $d=\dim(\mathcal V(N,M))$, the identity
$\dim(\mathcal V(N,M))=\dim(\mathcal V(N+1,M))$ shows
that the finite sequence $B[0],\dots,B[N+M+1]$ can be 
completed to a uniquely defined element $B\in\mathcal R$
which is of complexity $d$ and of saturation level $\leq N$. 
Use the natural isomorphisms between $\mathcal V(N,M),\mathcal V(N+1,M)$
and $\mathcal V(N,M+1)$ and the inclusions $\rho(s,t)\mathcal V(N,M+1)
\subset \mathcal V(N,M),\ 0\leq s,t\leq 1$ (where $\rho(s,t)$
acts in the obvious way) for writing
down a recursive presentation of $B$.

Check if $AB=1$ using the algorithms of Section \ref{sectsatlev}.

If yes, print the presentation found for $B$ and stop.

Otherwise, increase $N$ by $1$ and iterate the Loop. 
\newline
\noindent {\bf End of Loop}

\begin{rem} The above algorithm can perhaps be improved.
In particular, it is probably not necessary to consider 
all $(4^{M+1}-1)/3$ words of $\mathcal X^{\leq M}\subset
\mathcal  X^*$ during the computation of the dimension of $\mathcal V(N+1,M)$.
\end{rem}

\section{Modulo $2$ and $2-$valuations}\label{sectmod2}

\noindent
{\bf Proof of Theorem \ref{thmmod2}}
The infinite
symmetric Pascal matrix $P$ with coefficients $\left({i+j\choose i}\pmod 2
\right)\in\{0,1\}$ for $0\leq i,j$ defines a convergent element
(still denoted) $P\in\mathcal A$. It follows from Lucas's formula
(see Section \ref{subsectmod2}) that $P$ is a recurrence matrix of 
complexity $1$
recursively presented by
$$P=1,\left(\begin{array}{cc}P&P\\P&\end{array}\right)$$
(zero-entries are omitted).
The recurrence matrix $P\in\mathcal R$ has an $LU$ decomposition in
$\mathcal R$ given by the equality $P=LDL^t$
with $L,D\in\mathcal R$ of complexity $1$ 
defined by the recursive presentations
$$L=1,\left(\begin{array}{cc}L\\L&L\end{array}\right)\hbox{ and }
D=1,\left(\begin{array}{cc}D\\&-D\end{array}\right)$$
where $L$ is lower triangular unipotent and $D$ is diagonal.
An easy analysis of the coefficients of the diagonal recurrence matrix
$D$ ends the proof.\hfill$\Box$

\begin{rem} The convergent lower triangular
recurrence matrix $L$ and the convergent diagonal recurrence matrix
$D$ correspond to the infinite limit-matrices (still denoted) 
$L,D$ with coefficients given by $L_{i,j}=
\left({i+j\choose i}\pmod 2\right)\in\{0,1\}$ (for $0\leq i,j$) 
and $D_{n,n}=(-1)^{\nu_0+\nu_1+\nu_2+\dots}=(-1)^{\mathop{ds}(n)}$
where $n=\sum_{j\geq 0} \nu_j2^j\geq 0$ is a binary integer.
\end{rem}

\begin{rem} Recurrence matrices of complexity $1$ 
are, up to scalars, of the form 
$1,M,M\otimes M,M\otimes M\otimes M,\dots$
where $M$ is a complex $2\times 2$ matrix. 

It follows that the matrices $L$ and $D$ (and thus also $P=LDL^t$)
involved in the proof of Theorem \ref{thmmod2}
are invertible in $\mathcal R$. The recurrence matrix $D$ is its
own inverse. The inverse $L^{-1}$ of $L$ is recursively presented by
$L^{-1}=1,\left(\begin{array}{cc}L^{-1}\\
-L^{-1}&L^{-1}\end{array}\right)$.
\end{rem}

\begin{rem} The spectrum of recurrence matrices of complexity $1$ is 
easy to compute: Given a square matrix $M$ of size $d\times d$
with characteristic polynomial $\prod_{j=1}^d(t-\lambda_j)$ we have
$$\prod_{(j_1,\dots,j_n)\in\{1,\dots,d\}^n}(t-\lambda_{j_1}\cdots
\lambda_{j_n})$$
for the characteristic polynomial of the iterated tensor power
$M^{\otimes^n}$.
\end{rem}

\noindent
{\bf Proof of Theorem \ref{thmCRAS}}
The infinite matrix $V$ with coefficients
$$V_{s,t}=i^{v_2({s+t\choose s})}=
i^{\mathop{ds}(s)+\mathop{ds}(t)-\mathop{ds}(s+t)}$$
gives rise to a convergent element $V\in\mathcal A$.
A bit of work using Kummer's formulae (see Section \ref{subsec2val})
shows that $V=V_1$ is recurrence matrix given by
the recursive presentation
$$\begin{array}{l}
\displaystyle
V_1=1,\left(\begin{array}{cc}V_1&V_2\\
V_2&iV_1\end{array}\right)\\
\displaystyle
V_2=1,\left(\begin{array}{cc}V_1&-iV_1+(1+i)V_2\\
-iV_1+(1+i)V_2&-V_1\end{array}\right)\end{array}$$

We have $V=LDL^t\in\mathcal R$ 
with 
$L=L_1\in\mathcal R$ recursively presented by
$$\begin{array}{l}
\displaystyle L_1=1,\left(\begin{array}{cc}
L_1\\
L_3&L_4\end{array}\right)\\
\displaystyle L_2=0,\left(\begin{array}{cc}
0&-iL_2\\
-L_1+L_3&-iL_2-iL_4\end{array}\right)\\
\displaystyle L_3=1,\left(\begin{array}{cc}
L_1&L_2\\
-iL_1+(1+i)L_3&L_2+(1+i)L_4\end{array}\right)\\
\displaystyle L_4=1,\left(\begin{array}{cc}
L_1\\
(1-i)L_1+iL_3&L_4\end{array}\right)
\end{array}$$
and with diagonal $D=D_1\in\mathcal R$ recursively presented by
$$\begin{array}{l}
\displaystyle D_1=1,\left(\begin{array}{cc}
D_1\\
&D_2\end{array}\right)\\
\displaystyle D_2=-1+i,\left(\begin{array}{cc}
D_3&\\
&2D_1-D_2+2D_3\end{array}\right)\\
\displaystyle D_3=-1+i,\left(\begin{array}{cc}
D_3&\\
&-D_2\end{array}\right)\\
\end{array}$$
An analysis (left to the reader)
of the diagonal entries of $D_1$ ends the proof.\hfill$\Box$

\begin{rem} The recurrence matrices $L,D$ and $V=LDL^t$ are 
invertible in $\mathcal R$, see \cite{BacherCRAS}. 
\end{rem}

\section{Beeblebrox reduction}\label{sectBeeblebrox}

This section is devoted to proofs and complements
involving the Beeblerox reduction $\chi_B({n\choose k})$
of binomial coefficients.

\noindent
{\bf Proof of Theorem \ref{thmBeeblebroxformula}}
We have
$${2n\choose 2k}=\frac{(2n)\cdots(2n-2k+1)}{(2k)\cdots 1}=
{n\choose k}\frac{(2n-1)(2n-3)\cdots (2n-2k+1)}{
(2k-1)(2k-3)\cdots 1}$$
where both the numerator and the denominator of the fraction
$$F=\frac{(2n-1)(2n-3)\cdots (2n-2k+1)}{(2k-1)(2k-3)\cdots 1}$$
contain $k$ terms.
If $k$ is even, we have $F\equiv 1\pmod 4$
since the numerator and denominator of the fraction $F$
contain both $k/2$ factors $\equiv 1\pmod 4$ and $k/2$ factors $\equiv -1
\pmod 4$. If $k$ and $n$ are both odd,
the numerator and denominator of $F$ contain both 
$(k+1)/2$ factors $\equiv 1\pmod 4$ and $(k-1)/2$ factors
$\equiv -1\pmod 4$ and we have again $F\equiv 1\pmod 4$.
If $k$ is odd and $n$ is
even, then both binomial coefficients ${2n\choose 2k}$ and 
${n\choose k}$ are even and we have
thus $\chi_B ({2n\choose 2k})=\chi_B ({n\choose k})=0$.
This proves the first equality.

The binomial coefficient ${2n\choose 2k+1}$ 
is obviously even and this implies the second equality.

In the next case we have
$${2n+1\choose 2k}={n\choose k}\frac{(2n+1)(2n-1)\cdots (2n-2k+3)}{(2k-1)
(2k-3)\cdots (1)}$$
and the last fraction equals $1\pmod 4$ if $k$ is even. 

For $k$ odd and $n$ even, we have $\chi_B ({2n+1\choose 2k})=
\chi_B ({n\choose k})=0$ since 
${2n+1\choose 2k}\equiv {n\choose k}\equiv 0\pmod 2$.

For $n,k$ both odd, the correction $(-1)^k=-1$ equals the 
fraction  modulo $4$. This ends the proof of the third equality.

In the case of the last equality, we have 
$${2n+1\choose 2k+1}=\frac{2n+1}{2k+1}{2n\choose 2k}$$
which is even if $n\equiv 0\pmod 2$ and $k\equiv 1\pmod 2$.
If $n\equiv k\pmod 2$ then 
$\frac{2n+1}{2k+1}\equiv 1\pmod 4$. For $n$ odd and $k$ even 
we have $\frac{2n+1}{2k+1}\equiv -1=(-1)^{n(k+1)}\pmod 4$. 
The first equality and these observations complete the proof.
\hfill$\Box$

\subsection{Proof of Theorem \ref{thmBeeblebroxdet}}

{\bf Proof} As in Section \ref{subsecLU}, we consider the element (still 
denoted) $Z\in\mathcal A$ associated to the infinite matrix $Z$ with 
coefficients $Z_{s,t}=\chi_B({s+t\choose s}),\ 0\leq s,t$.

We have to show that $Z$ is a recurrence matrix and we have
to find a recursive presentation for $Z$. This can be done
in the following way: We consider {\it left shift-maps} 
$\lambda(0,0),\lambda(0,1),\lambda(1,0),\lambda(1,1)$ which associate
to $A=(A[0],A[1],\dots)\in\mathcal A$ the 
element $\lambda(s,t)A\in\mathcal A$
where $(\lambda(s,t)A)[n]$ is the submatrix of $A[n+1]$ 
corresponding to row-indices
$\equiv s\pmod 2$ and column-indices $\equiv
t\pmod 2$. A subspace $\mathcal V\subset \mathcal A$ is
{\it left-recursively closed} if it is invariant under all four left 
shift-maps and the {\it left-recursive closure} ${\overline A}^\lambda
\subset\mathcal A$ of $A\in\mathcal A$ is the smallest left-recursively
closed subspace containing $A$.
Theorem \ref{thmBeeblebroxformula}
implies that ${\overline Z}^\lambda$ is finite-dimensional.
One shows that $\mathop{dim}({\overline A}^\lambda)=
\mathop{dim}(\overline A)$ for all $A\in\mathcal A$ 
(see for example \cite{BacherRM} for the details).
This proves that $Z$ is a recurrence matrix. 
A little bit of work based on properties of the saturation 
index shows now that 
$Z=Z_1$ is given by the recursive presentation
$$\begin{array}{l}
\displaystyle Z_1=1,\left(\begin{array}{cc}Z_1&Z_2\\Z_3&0\end{array}\right)\\
\displaystyle Z_2=1,\left(\begin{array}{cc}Z_1&Z_2\\-Z_3&0\end{array}\right)\\
\displaystyle Z_3=1,\left(\begin{array}{cc}Z_1&-Z_2\\Z_3&0\end{array}\right)
\end{array}$$

We have the identity
$Z=LDL^t$ with $L=L_1\in\mathcal R$ lower triangular unipotent 
given by the recursive
presentation 
$$\begin{array}{ll}
\displaystyle L_1=1,\left(\begin{array}{cc}L_1&\\L_3&L_4\end{array}\right)
\qquad&
\displaystyle L_2=2,\left(\begin{array}{rr}-\frac{2}{3}L_1&2L_2\\
\frac{2}{3}L_3&2L_4\end{array}\right)\\
\displaystyle L_3=1,\left(\begin{array}{cc}L_1&L_2\\L_3&L_4\end{array}
\right)&
\displaystyle 
L_4=1,\left(\begin{array}{rr}L_1&\\\frac{1}{3}L_3&L_4\end{array}\right)
\end{array}$$

The diagonal matrix $D=D_1\in\mathcal R$ has recursive presentation
$$D_1=1,\left(\begin{array}{cc}D_1\\&D_2\end{array}\right),
D_2=-1,\left(\begin{array}{cc}3D_1\\&\frac{1}{3}D_2\end{array}\right)$$
An easy inspection of the diagonal entries of $D_1$
completes the proof.\hfill $\Box$

\begin{rem} The
birecursively closed subspaces of $\mathcal A$ appearing 
in Remark \ref{rembirec} are the recursively closed subspaces 
of $\mathcal A$ which are also left-recursively closed,
ie. invariant under all four left shift-maps mentionned above.
\end{rem}

\subsubsection{The group $\Gamma_L$}

All four recurrence matrices $L_1,L_2,L_3,L_4$ involved in our
recursive presentation of $L=L_1$ are invertible in $\mathcal R$.
Their inverses are given by 
$$L_1^{-1}=M_1,\ L_2^{-1}=-\frac{1}{2}M_3,\ L_3^{-1}=
-\frac{1}{2}M_2,\ L_4^{-1}=M_4$$
with $M_1,M_2,M_3,M_4$ recursively presented by
$$\begin{array}{ll}
\displaystyle M_1=1,\left(\begin{array}{cc}M_1&\\
M_3&M_4\end{array}\right)&
\displaystyle M_2=-2,\left(\begin{array}{cc}2M_1&2M_2\\
2M_3&2M_4\end{array}\right)\\
\displaystyle M_3=-1,\left(\begin{array}{rr}M_1&M_2\\
\frac{1}{3}M_3&-\frac{1}{3}M_4\end{array}\right)\qquad &
\displaystyle M_4=1,\left(\begin{array}{rr}M_1&\\
\frac{1}{3}M_3&M_4\end{array}\right)\end{array}$$

We have the curious inclusions
$$\begin{array}{ll}
\displaystyle \rho(0,0)L\in\mathbb QL_1,\quad &
\displaystyle \rho(0,1)L\in\mathbb QL_2,\\
\displaystyle \rho(1,0)L\in\mathbb QL_3,\quad &
\displaystyle \rho(1,1)L\in\mathbb QL_4\end{array}$$
for $L\in \overline{L_1}=\overline{L_2}=\overline{L_3}=
\overline{L_4}=\oplus_{j=1}^4 \mathbb C L_j$.
The analogous property holds also for their inverses
$$M_1=L_1^{-1},\ M_2=-2L_3^{-1},\ M_3=-2L_2^{-1},\ M_4=L_4^{-1}\ .$$

This suggests that it would perhaps be interesting to understand 
the group $\Gamma_L=\langle a,b,c,d\rangle\subset\mathcal R$
generated by the normalised recurrence matrices
$$a=L_1=1,\left(\begin{array}{cc}1\\1&1\end{array}\right),
\left(\begin{array}{cccc}1\\1&1\\1&2&1\\1&1&\frac{1}{3}&1\end{array}\right),
\dots$$
$$b=\frac{1}{2}L_2=1,
\left(\begin{array}{rr}-\frac{1}{3}&2\\\frac{1}{3}&1\end{array}\right),
\left(\begin{array}{rrrr}-\frac{1}{3}&0&-\frac{2}{3}&4
\\-\frac{1}{3}&-\frac{1}{3}&\frac{2}{3}&2\\\frac{1}{3}&\frac{2}{3}&1&0\\
\frac{1}{3}&\frac{1}{3}&\frac{1}{3}&1\end{array}\right),
\dots$$
$$c=L_3=1,
\left(\begin{array}{rr}1&2\\1&1\end{array}\right),
\left(\begin{array}{rrrr}1&0&-\frac{2}{3}&4
\\1&1&\frac{2}{3}&2\\1&2&1&0\\
1&1&\frac{1}{3}&1\end{array}\right),
\dots$$
$$d=L_4=1,
\left(\begin{array}{rr}1\\\frac{1}{3}&1\end{array}\right),
\left(\begin{array}{rrrr}1\\1&1\\
\frac{1}{3}&\frac{2}{3}&1\\
\frac{1}{3}&\frac{1}{3}&\frac{1}{3}&1\end{array}\right),
\dots$$
corresponding to the invertible ``projective'' elements $\mathbb Q^*L_i$.
In particular, it would be interesting to understand if 
$\Gamma_L$ is a linear group. This would certainly be implied by
the existence (which I ignore) of a natural integer
$N$ for which the projection $\Gamma_L\longrightarrow \pi_{\leq N}
(\Gamma_L)$ is one-to-one. 

L. Bartholdi communicated to me the following list implying 
all relations of length $\leq 12$ in $\Gamma_L$:
$$\begin{array}{l}
\displaystyle b d^{-1} c a^{-1},\\
\displaystyle  c a^{-1} c a^{-1},\\
\displaystyle  c a d^{-1} a^{-1} d^2 a^{-1} b^{-1},\\
\displaystyle  c a d^{-1} c^{-1} b d a^{-1} b^{-1},\\
\displaystyle  c d a^{-2} d a d^{-1} b^{-1},\\
\displaystyle  c a^2 d^{-1} a^{-2} d a d a^{-2} b^{-1},\\
\displaystyle  d^2 a d^{-2} b^{-1} c a d a^{-3},\\
\displaystyle  c d a d^{-2} a^{-1} d a d a^{-2} b^{-1}.\end{array}$$
He observed that they are all of the form $u_n v_n^{-1} w_n
x_n^{-1}$, where $u_n,\ v_n,\ w_n,\ x_n$ are positive words of length $n$
with respect to the generators $\{a,b,c,d\}$.
 
More generally, it should also be interesting to understand the subalgebra
$\mathcal L\subset \mathcal R$ generated by the recurrence matrices
$a^{\pm 1},b^{\pm 1},c^{\pm 1},d^{\pm 1}\in\mathcal R$.

The algebra $\mathcal L$ is of course a quotient of the group algebra 
$\mathbb C[\Gamma_L]$ and it would be interesting to describe the kernel
of the associated homomorphism.

The algebra $\mathcal L$ is a quotient of the free non-commutative 
algebra $$\mathcal N=
\mathbb C\langle A,A^{-1},B,B^{-1},C,C^{-1},D,D^{-1}\rangle$$
in eight free non-commutative variables $A^{\pm 1},B^{\pm 1},
C^{\pm 1},D^{\pm 1}$. 
The corresponding 
natural homomorphism $\pi:\mathcal N\longrightarrow \mathcal L$
(given by $Z^{\pm 1}\longmapsto z^{\pm 1}\in \mathcal R$ for $(Z,z)\in
\{(A,a),(B,b)(C,c),(D,d)\}$) 
factorises through the group algebra of the abstract group $\Gamma_L$.
The kernel $\mathcal I=\ker(\pi)\subset \mathcal N$ 
contains thus all relation of $\Gamma_L$ and in particular 
the trivial relations
$$AA^{-1}-1,\ BB^{-1}-1,\ CC^{-1}-1,\ DD^{-1}-1\ . $$

In order to gain some information on $\mathcal I$, we can
consider the morphism of algebras 
$\mu_1:\mathcal N\longrightarrow M_{2\times 2}(\mathcal N)$ given by
$$\begin{array}{ll}
\displaystyle 
A\longmapsto \left(\begin{array}{cc}A&0\\C&D\end{array}\right),& 
\displaystyle
A^{-1}\longmapsto \left(\begin{array}{cc}A^{-1}&0\\-B^{-1}&D^{-1}
\end{array}\right),\\
\displaystyle
B\longmapsto \left(\begin{array}{cc}-A/3&2B\\C/3&D\end{array}\right),&
\displaystyle
B^{-1}\longmapsto \left(\begin{array}{cc}-A^{-1}&2C^{-1}\\B^{-1}/3&D^{-1}/3
\end{array}\right),\\
\displaystyle
C\longmapsto \left(\begin{array}{cc}A&2B\\C&D\end{array}\right),&
\displaystyle
C^{-1}\longmapsto \left(\begin{array}{cc}-A^{-1}&2C^{-1}\\B^{-1}&-D^{-1}
\end{array}\right),\\
\displaystyle
D\longmapsto \left(\begin{array}{cc}A&0\\C/3&D\end{array}\right),&
\displaystyle
D^{-1}\longmapsto \left(\begin{array}{cc}A^{-1}&0\\-B^{-1}/3&D^{-1}
\end{array}\right),\end{array}$$

This morphism factors through $\pi$ and induces a homomorphism
of algebras $\overline \mu_1:\mathcal L\longrightarrow
M_{2\times 2}(\mathcal L)$
which removes simply the first matrix $X[0]$ from an
element  $X[0],X[1],\dots\in\mathcal L$.
This is due to the definition of $\mu_1$ which 
corresponds to the maps
$\mathcal R\longrightarrow M_{2\times 2}
(\mathcal R)$ given by $$X\longmapsto \left(\begin{array}{cc}
\rho(0,0)X&\rho(0,1)X\\\rho(1,0)X&\rho(1,1)X\end{array}\right)$$
for $X\in\{a^{\pm 1},b^{\pm-1},c^{\pm 1},d^{\pm 1}\}$.
We have thus in particular $\mu_1(\mathcal I)\subset
M_{2\times 2}(\mathcal I)$. Application of $\mu_1$ to some known
element in $\mathcal I$ can sometimes be used 
for the discovery of new elements in $\mathcal I$:
The computation 
$$\mu_1(AA^{-1})=\left(\begin{array}{cc}A&\\C&D\end{array}\right)
\left(\begin{array}{cc}A^{-1}\\
-B^{-1}&D^{-1}\end{array}\right)=
\left(\begin{array}{cc}AA^{-1}\\CA^{-1}-DB^{-1}&DD^{-1}\end{array}\right)$$
shows that application of $\mu_1$ to the trivial relation 
$AA^{-1}-1\in \mathcal I$ implies the already known inclusions 
$AA^{-1}-1,DD^{-1}-1\in\mathcal I$ and produces the non-trivial 
relation $CA^{-1}-DB^{-1}\in \mathcal I$. Since elements of 
$\mathcal I$ involving only two monomials induce relations on the group
$\Gamma_L$, we get the relation $bd^{-1}ca^{-1}$ in $\Gamma_L$
which is the first relation in Bartholdi's list.

``Iterating'' the map $\mu_1$ produces homomorphisms
$\mu_n:\mathcal N\longrightarrow M_{2^n\times 2^n}(\mathcal N)$
with similar properties. In particular, we have $\mu_n(\mathcal I)
\subset M_{2^n\times 2^n}(\mathcal I)$. Indexing the coefficients 
of matrices in $M_{2^n\times 2^n}(\mathcal N)$ by elements 
$\mathbf X\in\{X_{0,0},X_{0,1},
X_{1,0},X_{1,1}\}^n$, we get linear maps $\mu_{\mathbf X}:\mathcal N
\longrightarrow\mathcal N$ by considering the coefficient 
corresponding to $\mathbf X$ in $\mu_n(\mathcal N)$.
The reader should be warned that the map 
$\mathbf X\longmapsto \mu_{\mathbf X}\in
\mathop{End}(\mathcal N)$ is not a morphism of monoids
from $\{X_{0,0},X_{0,1},X_{1,0},X_{1,1}\}^*$ into 
$\mathop{End}(\mathcal N)$. The monoid 
generated by all maps $\mu_{\mathbf X},\mathbf X\in\{X_{0,0},X_{0,1},
X_{1,0},X_{1,1}\}^*,$ preserves however the ideal $\mathcal I$.

Denoting by $\mathcal N(\mathbb Z)$ the subring of noncommutative polynomials
with integral coefficients, relations of the form $r_1=r_2$ in $\Gamma_L$
are in bijection with pairs of roots in the infinite-dimensional
Euclidean lattice (with respect to the orthonormal basis given by 
monomials) $\mathcal I\cap \mathcal N(\mathbb Z)$.

The peculiar form of all relations in Bartholdi's 
list is partially explained by the formulae for $\mu_1$.
They imply that the maps $\mu_{\mathbf X}$
preserve sign structures:
The vector space spanned by the orbit under the monoid generated by
the maps $\mu_{\mathbf X}$ of a relation of the form 
$u_nv^{-1}_n=x_nw_n^{-1}$ contains only relations of the same form.

It would be interesting to know if the ideal $\mathcal I$ is
finitely generated as an $\{\mu_{\mathbf X}\}_{\mathbf X\in
\{X_{0,0},X_{0,1},X_{1,0},X_{1,1}\}^*}-$module:
Otherwise stated, 
does $\mathcal I$ contain a finite subset 
$\mathcal G$ such that $\mathcal I$ is the smallest bilateral ideal 
which contains $\mathcal G$ and which is preserved by all maps
$\mu_{\mathbf X}$?

\begin{rem} The techniques used in this Section can of course be
applied to other subsets in $\mathcal R$. One needs a (preferably
finitely generated) algebra $\mathcal S$ (eg. the algebra generated by 
suitable elements of a subgroup in $\mathcal R$) 
such that the recursive closure of every element in $\mathcal S$
is spanned by elements of $\mathcal S$. The choice of a generating set
$\mathcal G$ allows to consider the free non-commutative algebra 
$\mathcal N$ on $\mathcal G$ which gives rise 
to the natural surjective homomorphism $\pi:\mathcal N\longrightarrow
\mathcal S$. Choosing lifts of the shift maps, one 
constructs homomorphisms
$\mu_{\mathbb N}:\mathcal N\longrightarrow M_{2^n\times 2^n}(\mathcal N)$
giving rise to the linear maps $\mu_{\mathbf X}\in\mathop{End}(\mathcal N)$
preserving the bilateral ideal
$\mathcal I=\ker(\pi)$. 

The maps $\mu_{\mathbf X}$ can be choosen in order to 
preserve the grading of 
$\mathcal N$ if the generating set $\mathcal G$ spans a recursively
closed subspace.
\end{rem}

\subsubsection{The inverses of $D_1,D_2$ and the group $\Gamma_Z$}

The inverses of the diagonal recurrence matrices $D_1,D_2$ 
(defined by the decomposition $Z=LD_1L^t$ and by $D_2=\rho(1,1)D_1$)
are $E_1=D_1^{-1},\ E_2=D_2^{-1}$
recursively presented by
$$E_1=1,\left(\begin{array}{cc}E_1\\&E_2\end{array}\right),
E_2=-1,\left(\begin{array}{cc}\frac{1}{3}E_1
\\&3E_2\end{array}\right)$$

The inverses of the matrices $Z_1,Z_2,Z_3$ are $U_1=Z_1^{-1},
U_3=Z_2^{-1},U_2=Z_3^{-1}$ recursively presented by
$$\begin{array}{l}
\displaystyle
U_1=1,\left(\begin{array}{cc}0&U_2\\U_3&U_1-U_2-U_3\end{array}\right)\ ,\\
\displaystyle
U_2=1,\left(\begin{array}{cc}0&U_2\\-U_3&-U_1+U_2+U_3\end{array}\right)\ ,\\
\displaystyle
U_3=1,\left(\begin{array}{cc}0&-U_2\\U_3&-U_1+U_2+U_3\end{array}\right)\ .
\end{array}$$
Since $Z_1=Z_1^t$ and $Z_2=Z_3^t$, 
the group $\Gamma_Z=\langle a,b,c\rangle\subset
\mathcal R$ generated by $a=Z_1,b=Z_2,c=Z_3$ has 
an involutive automorphism given by $a\mapsto a^{-1},\ 
b\mapsto c^{-1},\ c\mapsto b^{-1}$.
Two relations in $\Gamma_Z$ 
are $(ab^{-1})^2$ and $ab=ca$.

Using the relation $a^2=cb$ following from the computation
$$a^2=aab^{-1}b=aba^{-1}b=caa^{-1}b=cb\ ,$$
every element of $\Gamma_Z$ can be expressed as an element of
$a^\epsilon \langle b,c\rangle$ with $\epsilon\in\{0,1\}$.

Since $\det(\pi_2(a))=-1$ and $\det(\pi_2(b))=\det(\pi_2(c))=1$,
the subgroup generated by $b,c$ is of index $2$ in $\Gamma_Z$.

\subsection{The triangular Beeblebrox matrix}

We define the lower triangular Beeblebrox matrix 
as the infinite lower triangular matrix with 
coefficients $L_{s,t}=\chi_B({s\choose t}),0\leq s,t,$ given by
the Beeblebrox reduction of binomial coefficients.

One of the main results of \cite{Gr1} states that any fixed
row of $L$ contains either no coefficients $-1$ or
the same number (given by a power of $2$) of coefficients
$1$ and $-1$. This can of course also be deduced from Theorem 
\ref{thmBeeblebroxformula} 
or by computing $LJ$ where $J$ is the ``recurrence vector''
obtained by considering the sequence of column vectors
$$(1),\ (1,1)^t,\ (1,1,1,1)^t,\ (1,1,1,1,1,1,1,1)^t,\dots\ .$$ 

The triangular Beeblebrox matrix $L$ defines a recurrence matrix
(still denoted) $L=L_1\in\mathcal R$ recursively presented by
$$L_1=1,\left(\begin{array}{cc}L_1\\L_2&L_3\end{array}\right),
L_2=1,\left(\begin{array}{cc}L_1\\L_2&-L_3\end{array}\right),
L_3=1,\left(\begin{array}{cc}L_1\\-L_2&L_3\end{array}\right)\ .$$

\subsubsection{The recurrence matrices $L_i^{-1}$}
The lower triangular recurrence matrices 
$L_1,L_2,L_3$ defined above are invertible in $\mathcal R$ 
with inverse elements $M_1=L_1^{-1},M_2=L_2^{-1},M_3=L_3^{-1}$
recursively presented by
$$M_1=1, \left(\begin{array}{cc}M_1\\M&M_3\end{array}\right),
M_2=1,\left(\begin{array}{rc} M_1\\-M&-M_3\end{array}\right),
M_3=1,\left(\begin{array}{rc}M_1&0\\-M&M_3\end{array}\right)$$
where $M=M_1-M_2-M_3$.

\begin{prop} \label{propgrouplowtriangBeeb} The map
$$L_1\longmapsto\left(\begin{array}{ccc}1&0&2\\0&1&0\\0&0&1\end{array}\right),
L_2\longmapsto\left(\begin{array}{ccc}0&1&1\\1&0&1\\0&0&1\end{array}\right),
L_3\longmapsto\left(\begin{array}{ccc}1&0&0\\0&1&2\\0&0&1\end{array}\right)$$
(where the three matrices correspond to the three affine 
maps $(x,y)\mapsto(x+2,y),
(x,y)\mapsto(y+1,x+1),(x,y)\mapsto(x,y+2)$ of $\mathbb R^2$)
defines a faithful linear representation of the group 
$\langle L_1,L_2,L_3\rangle\subset\mathcal R$ generated by 
$L_1,L_2,L_3$.

Moreover, the group homomorphism $L_i\longmapsto \pi_2(L_i)$
is faithful on $\langle L_1,L_2,L_3\rangle$.
\end{prop}

\noindent
{\bf Proof} We check that $L_1$ and $L_3$ commute.
They generate thus an abelian subgroup $\Gamma_a$
which is easily seen to be free abelian
of rank $2$ by considering the $4\times 4$ matrices 
$\pi_2(L_1)$ and $\pi_2(L_3)$. Checking the relations
$$L_2^2=L_1L_3,\ L_2L_3=L_1L_2,\ L_2L_1=L_3L_2$$
shows that $\Gamma_a$ is of index $2$ in $\langle L_1,L_2,L_3\rangle$
and these relations define the affine group of Proposition 
\ref{propgrouplowtriangBeeb}.

Faithfulness of the homomorphism $L_i\longmapsto \pi_2(L_i)$
follows from the observation
that the subgroup generated by $\pi_2(L_1),\pi_2(L_3)$
is free abelian of rank $2$ and does not contain $\pi_2(L_2)$.
\hfill$\Box$

\section{On the Jacobi-Dirichlet character}\label{sectjacDir}

This section contains the most important data for proving 
Theorem \ref{thmJacDir}. We omit the somewhat lengthy details.

Tedious work proving formulae analogous to Theorem 
\ref{thmBeeblebroxformula} or general principles show that the matrix 
$J=J_1$ with coefficients $\chi_J({s+t\choose s}),\ 0\leq s,t$ is 
of complexity $9$ and has  
a recursive presentation given by 
$$\begin{array}{ll}
\displaystyle
J_1=1,\left(\begin{array}{cc}J_1&J_2\\J_2^t&\end{array}\right),\\
\displaystyle
J_2=1,\left(\begin{array}{cc}J_3&J_4\\J_5&\end{array}\right),\qquad
J_2^t=1,\left(\begin{array}{cc}J_3^t&J_5^t\\J_4^t&\end{array}\right),\\
\displaystyle
J_3=1,\left(\begin{array}{cc}J_1&J_2\\-J_2^t&\end{array}\right),\qquad
J_3^t=1,\left(\begin{array}{cc}J_1&-J_2\\J_2^t&\end{array}\right),\\
\displaystyle
J_4=1,\left(\begin{array}{cc}J_3&J_4\\-J_5&\end{array}\right),\qquad
J_4^t=1,\left(\begin{array}{cc}J_3^t&-J_5^t\\J_4^t&\end{array}\right),\\
\displaystyle
J_5=-1,\left(\begin{array}{cc}J_3^t&J_5^t\\-J_4^t&\end{array}\right),\qquad
J_5^t=-1,\left(\begin{array}{cc}J_3&-J_4\\J_5&\end{array}\right),
\end{array}$$
where $X^t=X[0]^t,X[1]^t,X[2]^t,\dots$ for $X\in\mathcal A$.

The matrix $J$ has an $J=LDL^t$ decomposition in $\mathcal R$
with $L$ of complexity $20$ and $D$ of complexity $4$.

The lower triangular recurrence matrix 
$L=L_1$ involved in the decomposition $J=LDL^t$ has the
recursive presentation:
$$\begin{array}{ll}
L_1=1,\left(\begin{array}{cc}L_1&0\\ L_2&L_3\end{array}\right),&
L_2=1,\left(\begin{array}{cc}L_4&L_5\\ L_6&L_7\end{array}\right),\\
L_3=1,\left(\begin{array}{cc}L_8&0\\ L_9&L_{10}\end{array}\right),&
L_4=1,\left(\begin{array}{cc}L_1&L_{11}\\ L_2&L_3\end{array}\right),\\
L_5=2,\left(\begin{array}{cc}L_{12}&L_{13}\\ 0&0\end{array}\right),&
L_6=1,\left(\begin{array}{cc}L_4&L_{14}\\ L_6&L_7\end{array}\right),\\
L_7=1,\left(\begin{array}{cc}L_8-L_{12}&L_{15}\\ L_9&L_{10}\end{array}\right),&
L_8=1,\left(\begin{array}{cc}L_1&0\\ 1/3L_2&L_3\end{array}\right),\\
L_9=1/3,\left(\begin{array}{cc}L_4&3L_5\\ 1/3L_6&L_7\end{array}\right),&
L_{10}=1,\left(\begin{array}{cc}L_8&0\\ 1/3L_9&L_{10}\end{array}\right),\\
L_{11}=2,\left(\begin{array}{cc}L_{16}&L_{13}\\ L_{17}&L_{18}
   \end{array}\right),&
L_{12}=4/3,\left(\begin{array}{cc}L_{19}&4L_5\\ 0&0\end{array}\right),\\
L_{13}=4,\left(\begin{array}{cc}2/3L_{12}&2L_{13}\\ 0&0\end{array}\right),&
L_{14}=0,\left(\begin{array}{cc}L_{12}+L_{16}&L_{20}\\ 
   L_{17}&L_{18}\end{array}\right),\\
L_{15}=2,\left(\begin{array}{cc}-2/3L_{12}+L_{16}&L_{13}\\ 
   1/3L_{17}&L_{18}\end{array}\right),&
L_{16}=-2/3,\left(\begin{array}{cc}-2/3L_1&2L_{11}\\ 
   2/3L_2&2L_3\end{array}\right),\end{array}$$
$$\begin{array}{l}
L_{17}=2/3,\left(\begin{array}{cc}-2/3L_4&-4L_5+2L_{14}\\ 
   2/3L_6&2L_7\end{array}\right),\\
L_{18}=2,\left(\begin{array}{cc}-2/3L_8-2/3L_{12}&2L_{15}\\ 
   2/3L_9&2L_{10} \end{array}\right),\\
L_{19}=8/3,\left(\begin{array}{cc}2L_{19}&8/3L_5\\ 0&0\end{array}\right),\\
L_{20}=-4,\left(\begin{array}{cc}4/3L_{12}-2/3L_{16}&-2L_{13}+2L_{20}\\ 
   2/3L_{17}&2L_{18}\end{array}\right).
\end{array}$$

The four matrices $L_5,L_{12},L_{13},L_{19}$ span a somewhat 
trivial four-dimensional
subalgebra consisting only of matrix sequences with zero coefficients
except for the first row.

Consideration of the images $\rho(s,t)\overline L$ yields the following 
decomposition 
of the recursively closed vector space 
$\overline L=\bigoplus_{j=1}^{20}\mathbb C L_j$:
$$\begin{array}{l}
\displaystyle \rho(0,0)\mathcal V=\mathbb CL_1\oplus 
\mathbb CL_4\oplus \mathbb CL_8\oplus \mathbb CL_{12}\oplus \mathbb CL_{16}\oplus \mathbb CL_{19}\ ,\\
\displaystyle \rho(0,1)\mathcal V=
\mathbb C L_5\oplus \mathbb C L_{11}\oplus \mathbb C L_{13}
\oplus \mathbb C L_{14}\oplus \mathbb C L_{15}\oplus \mathbb C L_{20}\ ,\\
\displaystyle \rho(1,0)\mathcal V=\mathbb C
L_2\oplus \mathbb C L_6\oplus \mathbb C L_9\oplus \mathbb C L_{17}\ ,\\
\displaystyle \rho(1,1)\mathcal V=\mathbb C
L_3\oplus \mathbb C L_7\oplus \mathbb C L_{10}\oplus \mathbb C L_{18}\ .
\end{array}$$

I ignore if the vector space $\overline L$ contains generators of 
interesting algebras or groups.

The diagonal recurrence matrix
$D=D_1$ involved in $J=LDL^t$ is of complexity $4$
with recursive presentation 
$$\begin{array}{ll}
\displaystyle D_1=1,\left(\begin{array}{cc}
D_1&\\&D_2\end{array}\right),\quad &
\displaystyle D_2=-1,\left(\begin{array}{cc}
D_3&\\&D_4\end{array}\right),\\
\displaystyle D_3=3,\left(\begin{array}{cc}
3D_1&\\&1/3D_2\end{array}\right),\quad &
\displaystyle D_4=-1/3,\left(\begin{array}{cc}
3D_3&\\&1/3D_4\end{array}\right).\end{array}$$

\section{$q-$binomials}\label{sectqbin}

\noindent
{\bf Proof of Theorem \ref{thmformqbinn}} The result holds for $b=0$ 
or for $a\leq b$.
An induction on $a+b$ ends the proof. 
It splits into the four following subcases:

If $a,b\not\equiv 0\pmod n$:
$$\begin{array}{lcl}
{a\choose b}_q&=&\omega^b{a-1\choose b}_\omega+{a-1\choose
  b-1}_\omega\\
&=&{\lfloor a/n\rfloor\choose \lfloor b/n\rfloor}\left(\omega^b
{a-1\pmod n\choose b\pmod n}_\omega+{a-1\pmod n\choose b-1\pmod
  n}_\omega\right)\\
&=&{\lfloor a/n\rfloor\choose \lfloor b/n\rfloor}
{a\pmod n\choose b\pmod n}_\omega\end{array}$$

If $a\equiv 0\pmod n,b\not\equiv0\pmod n$:
$$\begin{array}{lcl}
{a\choose b}_q&=&\omega^b{a-1\choose b}_\omega+{a-1\choose
  b-1}_\omega\\
&=&{a/n-1\choose \lfloor b/n\rfloor}
\left(\omega^b{n-1\choose b\pmod n}_\omega+{n-1\choose (b-1)\pmod
    n}_\omega
\right)\\
&=&{a/n-1\choose \lfloor b/n\rfloor}
{n\choose b\pmod n}_\omega\end{array}$$
and ${n\choose b\pmod n}_\omega=0$ since $b\not\equiv 0\pmod n$
implies that it is divisible by the 
$n-$th cyclotomic polynomial.

If $a\not\equiv 0,b\equiv 0\pmod n$:
$$\begin{array}{lcl}
{a\choose b}_q&=&\omega^b{a-1\choose b}_\omega+{a-1\choose
  b-1}_\omega\\
&=&{a-1\choose b}_\omega+{\lfloor a/n\rfloor\choose b/n-1}{a-1\pmod n\choose n-1}_\omega\\
&=&{\lfloor a/n\rfloor\choose b/n}{a-1\pmod n\choose 0}_\omega+0\\
&=&{\lfloor a/n\rfloor\choose b/n}{a\pmod n\choose 0}_\omega
\end{array}$$

If $a\equiv b\equiv 0\pmod n$:
$$\begin{array}{lcl}
{a\choose b}_q&=&\omega^b{a-1\choose b}_\omega+{a-1\choose
  b-1}_\omega\\
&=&{a/n-1\choose b/n}{n-1\choose
  0}_\omega+{a/n-1\choose b/n-1}{n-1\choose n-1}_\omega\\
&=&{a/n\choose b/n}{0\choose 0}_\omega\end{array}$$
\hfill $\Box$

\begin{rem} I thank C. Krattenthaler for pointing out that 
Theorem \ref{thmformqbinn} follows also from the classical $q-$binomial
identity
$$\prod_{j=0}^{n-1}(1-q^jt)=\sum_{k=0}^n(-t)^kq^{{k\choose 2}}{n\choose k}_q
\ .$$
Setting $n=a$ and $q=\omega$ for $\omega$ a primitive
$d-$th root of $1$ and equating the coefficients of $t^b$ on both sides
implies the result easily for odd $d$. 
The case of $d$ even requires also a sign analysis.

The above identity follows by induction on $n$ from the easy computation
$$\begin{array}{l}
\prod_{j=0}^n(1-q^jt)=(1-q^nt)\sum_{k=0}^n(-t)^kq^{k\choose 2}{n\choose k}_q\\
=1+\sum_{k=1}^{n+1}(-t)^kq^{k\choose 2}\left({n\choose n-k}_q+q^{n+1-k}
{n\choose n+1-k}_q\right)\\
=\sum_{k=0}^{n+1}(-t)^kq^{k\choose 2}{n+1\choose k}_q\ .\end{array}$$
\end{rem}

%{\bf Proof of Corollary \ref{corformqbinn}} Theorem \ref{thmformqbinn}
%implies that the polynomial
%$$P(q)={a\choose b}_q-{\lfloor a/n\rfloor\choose \lfloor
%  b/n\rfloor}_{q^n}\ 
%{a\pmod n\choose b\pmod n}_q\in\mathbb Z[q]$$
%vanishes for $q$ a primitive $n-$th root
%of 1. This shows that $P(q)$ is divisible  
%by the cyclotomic polynomial $C_n(q)$.\hfill $\Box$ 

Theorem \ref{thmformqbinn} implies in particular that for 
$\omega=e^{2i\pi k/n}$, the 
complex numbers ${a\choose b}_{\omega},\ a,b\in\mathbb N$
belong to the finite subset $\cup_{0\leq b\leq a<n}
\mathbb R_{\geq 0}{a\pmod n\choose b\pmod n}_\omega$
of real half-lines in $\mathbb C$.

For $n=2$, the matrix ${a\choose b}_{-1},\ 0\leq a,b<2$ is given by
$$\left(\begin{array}{cc}1&0\\1&1\end{array}\right)\ .$$
This implies that ${a\choose b}_{-1}\in\mathbb N$ for all $a,b\in\mathbb N$.

For $n=4$ we get 
${a\choose b}_i\in\mathbb N\cup i\mathbb N\cup(1+i)\mathbb N$
since we have
$$\left(\begin{array}{cccc}
1\\1&1\\
1&1+i&1\\
1&i&i&1\end{array}\right)$$
for the matrix with coefficients ${a\choose b}_i,\ 0\leq a,b<4$.

\subsection{Tensor products}\label{subsecttensor}

Replacing coefficients in the ground-field $\mathbb C$ by
multiples of a fixed matrix $X$ of size $2^K\times 2^K$,
one can consider the tensor product 
$$A\otimes X=(A[0]\otimes X,A[1]\otimes X,\dots)$$
of $A\in\mathcal A$ or $A\in\mathcal R$ with $X$. Such an element has an $LU$
decomposition $A\otimes X=(L'\otimes L_X)(U'\otimes U_X)$
involving elements of the same form if and only if we have decompositions 
$A=L'U'\in\mathcal R$ and $X=L_XU_X$.

\begin{rem} More generally, one can consider the quotient algebras
$\mathcal A/\mathcal{FS}$ and $\mathcal R/\mathcal{FS}$
where $\mathcal{FS}$ is the ideal of all matrix sequences 
in $\mathcal A$ which involve only finitely many non-zero matrices.
Elements of the quotient algebra $\mathcal R/\mathcal{FS}$
can be represented as linear combinations of suitable elements in 
$\cup_{K\in\mathbb N}\mathcal R\otimes M_{2^K\times 2^K}$
and such representations are sometimes simpler than 
recursive presentations of preimages in $\mathcal R$.
\end{rem}

\subsection{The Beeblebrox reduction of ${s+t\choose s}_{-1}$}
\label{subsectionbeebred}

We denote by $Z'$ the infinite symmetric matrix
with coefficients $Z'_{s,t}=\chi_B({s+t\choose s}_{-1})$
given by the Beeblebrox reduction of
$q-$binomials evaluated at $q=-1$.

Theorem \ref{thmformqbinn} and Section \ref{subsecttensor} imply that 
$Z'=L'D'(L')^t$ where 
$$Z'=Z\otimes \left(\begin{array}{cc}1&1\\1\end{array}\right),\ 
L'=L\otimes \left(\begin{array}{cc}1\\1&1\end{array}\right),\ 
D'=D\otimes \left(\begin{array}{cc}1\\&-1\end{array}\right)$$
(the tensor product $X\otimes M$ denotes the matrix(-sequence) obtained
by replacing a scalar entry $\lambda$ of $Z$ by the $2\times 2$ matrix $M$)
with $Z,L,D$ as in the proof of Theorem \ref{thmBeeblebroxdet}.

In particular, using Theorem \ref{thmBeeblebroxdet}
and the $LU$ decomposition of 
$\left(\begin{array}{cc}1&1\\1\end{array}\right)=\left(\chi_B({s+t
\choose s}_{-1})\right)_{0\leq s,t\leq 1}$,
one can easily write down a formula for
$\det(Z'(n))\in\pm3^{\mathbb N}$ with
$Z'(n)$ denoting the symmetric $n\times n$ 
submatrix consisting of the first $n$ rows and columns of $Z'$.

\begin{rem} The case of the matrix (with coefficients in $\{0,1\}$)
obtained by reducing ${s+t\choose s}_{-1}$ modulo $2$ yields nothing 
new since ${s+t\choose s}_{-1}\equiv 
{s+t\choose s}_1={s+t\choose s}\pmod 2$.
\end{rem}

\subsection{Reduction modulo $2$ and Beeblebrox reduction of 
${s+t\choose s}_i$}

Let $M'$ be the symmetric matrix with coefficients
$$\psi({s+t\choose s}_i)\in\{0,\pm 1,\pm x,\pm y\},\ 0\leq s,t$$ 
where
$$\psi(\xi)=\left\lbrace
\begin{array}{ll}
\gamma(\xi)&\hbox{if }\xi\in\mathbb N\\
\gamma(a)x&\hbox{if }\xi=ai\in i \mathbb N\\
\gamma(a)y&\hbox{if }\xi=a(1+i)\in (1+i)\mathbb N\end{array}\right.$$
where $\gamma:\mathbb N\longrightarrow \{0,\pm 1\}$ is either 
the reduction modulo $2$ with values in $\{0,1\}$ or the Beeblebrox
reduction $\chi_B$.

As in section \ref{subsectionbeebred} we have 
$M'=L'D'(L')^t$ where
$$\begin{array}{l}
\displaystyle M'=M\otimes \left(\begin{array}{cccc}1&1&1&1\\
1&y&x\\
1&x\\1\end{array}\right)\ ,\\
\displaystyle L'=L\otimes
\left(\begin{array}{cccc}
1\\
1&1\\
1&\frac{1-x}{1-y}&1\\
1&\frac{1}{1-y}&\frac{y-x}{x^2-2x+y}&1\end{array}\right)\ ,\\
\displaystyle 
D'=D\otimes \left(\begin{array}{cccc}
1\\
&y-1\\
&&\frac{x^2-2x+y}{1-y}\\
&&&\frac{-x^2}{x^2-2x+y}\end{array}\right)\ .\end{array}$$
$M,L,D$ are given by the matrices $P,L,D$, respectively $Z,L,D$ 
occuring in the proof of Theorem \ref{thmmod2}, respectively 
\ref{thmBeeblebroxdet}, if $\gamma$ is the reduction 
modulo $2$, respectively the Beeblebrox reduction.

It follows that the determinant $\det(M'(n))$ of 
the finite matrix $M'(n)$ consisting of the first $n$
rows and columns of $M'$ is of the form
$$\pm 3^{\mathbb N}x^{2\mathbb N}(y-1)^{\{0,1\}}(x^2-2x+y)^{
\{0,1\}}$$
with powers of $3$ only involved if $\gamma$ is the Beeblebrox reduction.
The factor $(y-1)$ appears if and only if $n\equiv 2\pmod 4$ and the
factor $(x^2-2x+y)$ appears if and only if $n\equiv 3\pmod 4$.
Using Theorem \ref{thmmod2}, respectively Theorem \ref{thmBeeblebroxdet},
it is easy to write down a formula for $\det(M'(n))$.

{\bf Acknowledgements} I thank L. Bartholdi,
M. Brion, P. de la Harpe and C. Krattenthaler for helpful
comments and discussions.

\noindent Roland BACHER

\noindent INSTITUT FOURIER

\noindent Laboratoire de Math\'ematiques

\noindent UMR 5582 (UJF-CNRS)

\noindent BP 74

\noindent 38402 St Martin d'H\`eres Cedex (France)
\medskip

\noindent e-mail: Roland.Bacher@ujf-grenoble.fr

\end{document}